\documentclass[12pt, a4paper]{amsart}
\usepackage{amsmath,amsfonts,amsthm,amssymb,bbm,xcolor,comment,epsfig,epstopdf,graphicx,eucal,natbib,dsfont,accents,url,mathrsfs,upgreek,hyperref,xr,caption}
\usepackage[top=3.5cm, bottom=3.5cm, left=3cm, right=3cm]{geometry}

\graphicspath{{Fig/}}

\newtheorem{thm}{Theorem}[section]
\newtheorem{lem}[thm]{Lemma}
\newtheorem{cor}[thm]{Corollary}
\numberwithin{figure}{section}
\numberwithin{table}{section}

\theoremstyle{remark}
\newtheorem{rem}[thm]{Remark}

\DeclareMathOperator{\cov}{Cov}

\newcommand*{\dt}[1]{%
  \accentset{\mbox{\large\bfseries .}}{#1}}
\newcommand{\dd}{\mathrm{d}}
\newcommand{\E}{\mathbb{E}}

\begin{document}
\title[Goodness-of-Fit Testing for Point Processes]{Goodness-of-Fit Testing for Point Processes in\\ Large Populations}
\author[S. U. Can]{Sami Umut Can}
\address{Department of Quantitative Economics, University of Amsterdam, P.O. Box 15867, 1001 NJ Amsterdam, The Netherlands}
\email{S.U.Can@uva.nl}
\author[E. V. Khmaladze]{Estate V. Khmaladze}
\address{School of Mathematics \& Statistics, Victoria University of Wellington, P.O. Box 600, Wellington, New Zealand}
\email{Estate.Khmaladze@vuw.ac.nz}
\author[R. J. A. Laeven]{Roger J.~A. Laeven}
\address{Department of Quantitative Economics, University of Amsterdam, P.O. Box 15867, 1001 NJ Amsterdam, The Netherlands}
\email{R.J.A.Laeven@uva.nl}

\begin{abstract}
Suppose we have an observed path from a point process counting event occurrences in a large population. Based on the observed path, we would like to test the null hypothesis that the conditional intensity of the point process belongs to a particular parametric family. We propose a novel approach to conducting such goodness-of-fit tests. The idea is to construct a unitary transformation of a natural parametric testing process such that it converges weakly to a ``standard'' target process, independent of the particular parametric form assumed under the null hypothesis. This transformation therefore paves the way for asymptotically distribution-free goodness-of-fit testing of parametric point processes. We demonstrate the good finite-sample performance of our approach through Monte Carlo simulations of Aalen-type survival processes, without and with censoring, mixture cure models, and software reliability models, and we illustrate its applicability with observed human lifetimes as well as real software failures.\\[2mm]
\noindent \textbf{\keywordsname:}\ Population point processes; conditional intensity; event history analysis; survival models; 
goodness-of-fit testing.\\[2mm]
\end{abstract}

\maketitle

\section{Introduction}\label{sec:intro}

There is a rich literature on the use of temporal point processes, or counting processes, in event history analysis as well as in the related fields of survival analysis and duration analysis, tracing back to the seminal work of \cite{aalen:1975}. A (simple) point process can be characterized by its conditional intensity process
describing the random arrival rate of events. The conditional intensity process in event history analysis, and in actuarial, demographic, and epidemiological statistics in particular, is often specified parametrically. A multitude of parametric models exist for the conditional intensity process. Therefore, a natural problem that arises is that of testing a parametric intensity model's goodness-of-fit to a given data sample, consisting of an observed path from the population point process.

In the---often fictitious---case where the parameter values are known rather than estimated, the well-known random \textit{time transformation}, which transforms the point process into a standard unit-rate Poisson process by a stochastic time change, provides a natural approach to analyzing an intensity model's goodness-of-fit; see e.g., \cite{meyer:1971}, \cite{bremaud:1972}, \cite{papangelou:1972}, \cite{aalen:hoem:1978} and \cite{brown:nair:1988}. Indeed, in the case of known parameters, and hence simple null hypotheses, goodness-of-fit tests can be based upon assessing the closeness of the time-transformed point process to a standard Poisson process. However, in the more realistic case where parameters are estimated, i.e., under a parametric null hypothesis, the stochastic time change is rendered invalid and this approach generally breaks down.

In this paper, we propose a novel approach to conducting goodness-of-fit tests for parametric models of the conditional intensity process. The idea is to consider a natural testing process based on the difference between the observed point process and its integrated intensity, which is estimated parametrically. 
We then describe a \emph{unitary transformation} that, as we establish, turns this process into another one that has a ``standard'' weak limit (as the population size grows) when the parametric null model is correct.\footnote{The \emph{unitary} feature of the transformation implies that it is an inner product-preserving bijection (or isometric isomorphism); see Section \ref{sec:uni} for details.} Because this weak limit is a standard process, independent of the parametric family that is tested for,
and of that model's true parameter vector, any testing functional of this process needs to be tabulated only once for use in all such goodness-of-fit testing problems. In other words, our approach is asymptotically distribution-free. We do not propose one specific test statistic; instead we construct a full test \emph{process} that always has the same weak limit under the null hypothesis and can therefore serve to generate a plenitude of asymptotically distribution-free test statistics.

We illustrate the applicability of our approach in Monte Carlo simulations for \cite{aalen:1978}-type survival processes, without and with censoring, mixture cure models, and software reliability models. The simulation results demonstrate that statistical tests based on our approach have good size and power properties in a wide variety of realistic settings.

Finally, we apply our testing procedure to real-life data on human lifetimes and software failures. When applied to mortality data of Luxembourgish males, our testing procedure rejects the Aalen survival model with Weibull force of mortality, but not the one with Gompertz force of mortality. When applied to software failure data reported by the Centre for Software Reliability, our approach rejects the classical Jelinski-Moranda model of software reliability, in favor of the more flexible Littlewood model.

Several approaches to goodness-of-fit testing of parametric intensity models for point processes exist. One can use specialized bootstrap procedures for point-process based statistics in the context of event history and survival analysis, to obtain approximate $p$-values for test statistics under the null hypothesis; see, e.g., \cite{lin:wei:ying:1993}, the monograph \cite{martinussen:scheike:2006}, \cite{dietrich:dobler:degunst:2025} and the many references therein. These resampling methods are quite flexible, but do not share the asymptotically distribution-free feature of our approach. Alternatively, one may use the transformation introduced in \cite{khmaladze:1981}, to obtain an empirical process that converges to a standard Brownian motion under the null hypothesis, thus enabling distribution-free testing; see also the discussion in Section~\ref{sec:point}. Our approach avoids resampling procedures and is arguably simpler and more flexible than approaches based on \cite{khmaladze:1981}. 

This paper is organized as follows. Section~\ref{sec:point} considers the (compensated) point process under a parametric null hypothesis. We will see that its weak limit depends on the parametric family that is tested for as well as on the true parameter values, hence is not suitable for distribution-free testing. Section~\ref{sec:uni} describes a unitary transformation that turns the limiting process into a standard target process and plays a pivotal role in our testing procedure. Section~\ref{sec:gof} applies an empirical version of this transformation and establishes its weak convergence under the null hypothesis. In Section~\ref{sec:mc}, we analyze the finite-sample performance of our testing procedure. Section~\ref{sec:data} applies our approach to two real-world data sets. Proofs, technical details, and additional simulation results are provided in four appendices following the main text. {\tt{R}} code to implement the procedure developed in this paper is available from the authors upon request.

\section{Point Processes, Martingales, and Parametric Estimation}\label{sec:point}

Suppose we are given an observed path from a (simple) temporal point process $N_{n}$ on a time interval $[0,T]$. We interpret $N_{n}(t)$ as the number of ``events'' (such as deaths or failures) that have occurred by time $t$ in a population of initial size $n$, with $N_{n}(0)=0$. We denote by $\lambda_{n}(t)$ the conditional intensity process associated with $N_{n}(t)$. It is the local conditional expectation of the increments of $N_{n}$ given its strict history. Formally, 
\begin{equation}
\lambda_{n}(t)=\lim_{\Delta\downarrow 0}\frac{\mathbb{E}\left[N_{n}(t+\Delta)-N_{n}(t)\vert\mathcal{F}_{t-}\right]}{\Delta},
\label{eq:lam}
\end{equation}
where $\mathcal{F}_{t-}$ represents the strict past of $N_{n}(t)$,
that is, $\mathcal{F}_{t-}=\sigma\left(N_{n}(s): 0\leq s< t\right)$. Furthermore, we denote by $\Lambda_{n}(t)$ the compensator associated with $N_{t}$, which takes the form of the integrated conditional intensity, that is:
\begin{equation}
\Lambda_{n}(t)=\int_{0}^{t}\lambda_{n}(s)\,\mathrm{d}s.
\label{eq:Lam}
\end{equation}
The behavior of the compensator is usually smoother than that of the associated point process. As is well-known, the difference $N_n(t) - \Lambda_n(t)$ is a zero-mean martingale and is termed the compensated point process. We refer to \cite{karr:1991} and \cite{daley:verejones:2003} for textbook treatments
of these and other concepts related to point processes.

In the literature, a variety of parametric models for the conditional intensity process exists. Suppose we select a parametric family $\{\lambda_{n,\theta}:\theta\in\Theta\}$, where $\theta$ is an $m$-dimensional parameter vector and $\Theta$ is an open subset of $\mathbb{R}^{m}$. (The openness of $\Theta$ implies that each $\theta$ has a neighborhood in $\Theta$, so that differentiation at $\theta$ will not induce any boundary effects.)
Examples of such parametric families include
\begin{equation*}
\lambda_{n,\theta}(t)=\frac{f_{\theta}(t)}{1-F_{\theta}(t)}Y_n(t),\qquad t\geq 0,  
\end{equation*}
where $f_{\theta}$ and $F_{\theta}$ are the probability density function (pdf) and cumulative distribution function (cdf)
of some lifetime or waiting time distribution on $[0,\infty)$, and $Y_n(t)$ is the number of individuals ``at risk'' for event occurrences just before time $t$. Given the single observed path from $N_{n}(t), \, t \in [0,T]$, our aim is then to decide whether the parametric model for $\lambda_{n}$ fits the data well. That is, we would like to test the null hypothesis
\begin{equation}
\mathcal{H}_{0}:\lambda_{n}\in\mathcal{L}_{n}=\{\lambda_{n,\theta}:\theta\in\Theta\},
\label{eq:H0}
\end{equation}
against the alternative hypothesis $\mathcal{H}_{1}: \lambda_{n}\notin\mathcal{L}_{n}=\{\lambda_{n,\theta}: \theta\in\Theta\}$. 
If the null hypothesis holds, then there is a $\theta_{0}\in\Theta$
such that $\lambda_{n}=\lambda_{n,\theta_{0}}$.

Let us, for illustration purposes, first consider the much simpler problem of testing the goodness-of-fit of a fully specified conditional intensity process rather than of a parametric family.
That is, consider the simple null hypothesis
\begin{equation}
\mathcal{H}_{0}:\lambda_{n}=\lambda_{n}^{*},
\label{eq:H0simple}
\end{equation}
where $\lambda_{n}^{*}$ has no free parameters. Two broad approaches to this problem exist: time transformation and martingale transformation. We consider the latter approach as it is the most relevant for our purposes; the former approach is briefly described and referenced in the Introduction.
The well-known martingale central limit theorem (e.g., Ch.~5 of \citealp{fleming:harrington:1991}) implies that for a point process $N_{n}(t)$ with conditional intensity $\lambda_{n}(t)$,
in regular cases as $n\rightarrow\infty$, the weak limit of the process
\begin{align}
W_{n}(t)=\frac{1}{\sqrt{n}}\left[N_{n}(t)-\Lambda_{n}(t)\right] = \frac{1}{\sqrt{n}}\left[N_{n}(t)-\int_{0}^{t}\lambda_{n}(s)\,\mathrm{d}s\right],\quad t\geq 0,
\label{eq:Wnfullyspecified}
\end{align}
is a Brownian motion $W_{B}$ with ``time'' $B(t)=\lim_{n\rightarrow\infty}\frac{1}{n}\Lambda_{n}(t)$, that is, a zero-mean Gaussian process $W_B$ on $[0,\infty)$ with covariance $\cov[W_B(s), W_B(t)] = B(s \wedge t)$. 
Therefore, goodness-of-fit tests for the simple null hypothesis \eqref{eq:H0simple} can be performed by constructing the path
\begin{equation}
W_{n}^{*}(t)=\frac{1}{\sqrt{n}}\left[N_{n}(t)-\Lambda_{n}^{*}(t)\right],\qquad t \in [0,T],
\label{eq:Wnsimple}
\end{equation}
and assessing whether this path is consistent with the statistical behavior of $W_{B^*}$, where
$B^{*}(t)=\lim_{n\rightarrow\infty}\frac{1}{n}\Lambda_{n}^{*}(t)$.
This can be done, for example, by assessing whether the normalized path
\begin{equation*}
W_{0,n}^{*}(t)=\int_{0}^{t}\frac{1}{\sqrt{\beta^{*}(s)}}\,\mathrm{d}W_{n}^{*}(s),\qquad t \in [0,T],
\end{equation*}
with $\beta^{*}(t)=\lim_{n\rightarrow\infty}\frac{1}{n}\lambda_{n}^{*}(t)$,
is consistent with the behavior of a standard Brownian motion.

Now let us turn to the parametric null hypothesis \eqref{eq:H0}.
Given a sample from $N_{n}(t)$ on $0\leq t\leq T$, we can mimic the approach above and consider the compensated process
\begin{equation}\label{eq:Wnparametric}
\widehat{W}_{n}(t) =\frac{1}{\sqrt{n}}\left[N_{n}(t)-\Lambda_{n,\widehat{\theta}}(t)\right] = \frac{1}{\sqrt{n}}\left[N_{n}(t)-\int_{0}^{t}\lambda_{n,\widehat{\theta}}(s)\,\mathrm{d}s\right], \quad t \in [0,T],
\end{equation}
analogous to \eqref{eq:Wnsimple}, where $\widehat{\theta}$ denotes the maximum likelihood estimator (MLE) of $\theta$, maximizing the log-likelihood function
\begin{equation}\label{MLE}
L(\theta) = \int_0^T \log \lambda_{n,\theta}(t) \, \dd N_n(t) - \int_0^T \lambda_{n, \theta}(t) \, \dd t.
\end{equation}
In order to use the process $\widehat{W}_n$ for goodness-of-fit testing, we first need to establish its asymptotic behavior as $n \to \infty$ under the null hypothesis \eqref{eq:H0}. Thankfully, this asymptotic behavior is well studied in the literature; see Ch.\ VI of \cite{andersen:etal:1993} and \cite{sun:etal:2001}. 
We will first state some regularity assumptions on the family of conditional intensities $\mathcal{L}_n = \{\lambda_{n,\theta}: \theta \in \Theta\}$, and then state the asymptotic behavior of $\widehat{W}_n$ under those assumptions, and under the null hypothesis \eqref{eq:H0}.

\bigskip

\textbf{A1.} There exists a neighborhood $\Theta_0$ of $\theta_0$ such that for all $\theta \in \Theta_0$, $t \in [0,T]$, and for all $i,j,k \in \{1,\ldots,m\}$, the partial derivatives
\[
\frac{\partial}{\partial \theta_i} \lambda_{n,\theta}(t), \quad \frac{\partial^2}{\partial \theta_i \, \partial \theta_j} \lambda_{n,\theta}(t), \quad \frac{\partial^3}{\partial\theta_i\,\partial \theta_j \, \partial \theta_k} \lambda_{n,\theta}(t),
\]
as well as
\[
\frac{\partial}{\partial \theta_i} \log\lambda_{n,\theta}(t), \quad \frac{\partial^2}{\partial \theta_i \, \partial \theta_j} \log\lambda_{n,\theta}(t), \quad \frac{\partial^3}{\partial\theta_i\,\partial \theta_j \, \partial \theta_k} \log\lambda_{n,\theta}(t),
\]
exist and are continuous in $\theta$. 
Furthermore, these partial derivatives are bounded in probability on $[0,T] \times \Theta_0$.

\bigskip

\textbf{A2.} For each $\theta \in \Theta_0$, there exist functions $\alpha_\theta: [0,T] \to \mathbb{R}^m$ and $\beta_\theta: [0,T] \to [0,\infty)$ such that
\begin{equation}\label{eq:alpha-beta}
\frac{\dt{\lambda}_{n,\theta}(t)}{\lambda_{n,\theta}(t)} \stackrel{P}{\to} \alpha_{\theta}(t),
\qquad
\frac{\lambda_{n,\theta}(t)}{n} \stackrel{P}{\to} \beta_{\theta}(t),
\end{equation}
uniformly in $t \in [0,T]$, with
\[
\dt{\lambda}_{n,\theta}(t)
:=\left(\frac{\partial \lambda_{n,\theta}(t)}{\partial \theta_{1}},\ldots,\frac{\partial \lambda_{n,\theta}(t)}{\partial \theta_{m}}\right)^{\top}.
\]

\medskip

\textbf{A3.} For each $\theta \in \Theta_0$, the matrix $\Gamma_\theta := \int_0^T \alpha_\theta(t) \, \alpha^\top_\theta(t) \, \beta_\theta(t) \, \dd t$ is positive definite.

\bigskip

Note that our Assumptions A1--A3 imply Condition VI.1.1 of \cite{andersen:etal:1993}, and therefore we have $\widehat\theta \stackrel{P}{\to} \theta_0$ by Theorem VI.1.1 there. Moreover, we obtain the following theorem regarding the asymptotic behavior of the process $\widehat{W}_n$ in \eqref{eq:Wnparametric}, as a direct consequence of Theorem 1 in \cite{sun:etal:2001}; also see Section VI.3.3.1 of  \cite{andersen:etal:1993} for a very similar result.

\begin{thm}\label{thm:1}
Under the Assumptions A1--A3 stated above, and under the null hypothesis \eqref{eq:H0}, the process $\widehat{W}_{n}$ defined in \eqref{eq:Wnparametric} converges weakly to
\begin{equation}
V_{B_{\theta_{0}}}(t)
:=W_{B_{\theta_{0}}}(t)
-\int_{0}^{t}\alpha^{\top}_{\theta_{0}}(s)\beta_{\theta_{0}}(s)\,\mathrm{d}s\,\Gamma_{\theta_{0}}^{-1}\int_{0}^{T}\alpha_{\theta_{0}}(s')\,\mathrm{d}W_{B_{\theta_{0}}}(s'),
\label{eq:Wparametric}
\end{equation}
in the Skorohod space $D([0,T])$,
where $B_{\theta_{0}}(t)=\int_0^t \beta_{\theta_0}(s)\,\dd s$.
\end{thm}

\begin{rem}
Theorem 1 of \cite{sun:etal:2001} considers multivariate point processes. In the univariate case, the process $n^{-1/2}\mathbf{M}(t; \hat{\boldsymbol{\uptheta}})$ defined there reduces to our process $\widehat{W}_n$, whereas the process $\boldsymbol{\upeta}(t; \boldsymbol{\uptheta}_0)$ defined there is the same as our process $V_{B_{\theta_{0}}}(t)$, with $f_1(t; \boldsymbol{\uptheta}_0)$ corresponding to our $\beta_{\theta_0}(t)$, $\mathbf{f}'_1(t;\boldsymbol{\uptheta}_0)/f_1(t;\boldsymbol{\uptheta}_0)$ corresponding to our $\alpha_{\theta_0}(t)$, and the matrix $\mathscr{T}$ corresponding to our $\Gamma_{\theta_0}$.
\end{rem}

It is clearly visible from Theorem~\ref{thm:1} that the behavior of the limiting process $V_{B_{\theta_{0}}}$ depends not only on the parametric family $\mathcal{L}_{n}=\{\lambda_{n,\theta}: \theta\in\Theta\}$ that is tested for, but also on the unknown true value of the parameter vector $\theta_{0}$. 
Therefore, the process $\widehat{W}_{n}$ and its weak limit $V_{B_{\theta_{0}}}$ under the null hypothesis do not provide a suitable basis for distribution-free goodness-of-fit testing. 
One way to solve this problem would be to apply the innovation martingale transformation described in \cite{khmaladze:1981,khmaladze:1988,khmaladze:1993}, which allows us to transform the limiting process \eqref{eq:Wparametric} into a standard Brownian motion on $[0,T]$. 
Then the same transformation applied to $\widehat{W}_n$ will be asymptotically a standard Brownian motion under the null hypothesis, which leads to asymptotically distribution-free testing. This idea has been explored in the literature in different contexts; see, e.g., Section VI.3.3.4 of \cite{andersen:etal:1993}, \cite{sun:etal:2001}, and \cite{bcl:2025}.

In the present work, we develop a novel approach, inspired by \cite{khmaladze:2016}, to transform $V_{B_{\theta_{0}}}$ in \eqref{eq:Wparametric} into a different ``standard'' target process. \cite{khmaladze:2016} proposes unitary transformations for Brownian bridges, in the context of testing multivariate distribution functions. Our approach relies on constructing a suitable unitary transformation of $V_{B_{\theta_{0}}}$ and of its empirical counterpart $\widehat{W}_n$. 
Under our approach, the target process will not be a standard Brownian motion. 
Instead, it will simply be another process of the form \eqref{eq:Wparametric}, one which we will define as the ``standard'' among all processes of this form. 
While this choice can, in principle, be arbitrary, a natural candidate is the process that arises as the limit of $\widehat{W}_n$ in \eqref{eq:Wnparametric} when $N_n$ is a standard Poisson process embedded in a convenient parametric family of conditional intensities. 
This will indeed be our choice for the target process. 
Since this target process is independent of the particular testing problem we are considering, and since we can explicitly calculate distributions of functionals from this target process, we will have a suitable basis for constructing asymptotically distribution-free goodness-of-fit tests, as we will see.   

\section{A Unitary Transformation}\label{sec:uni}

\subsection{Preliminaries}
We imagine two point processes $N_{\lambda_n}(t)$ and $N_{\mu_n}(t)$ for $t \in [0,T]$, with respective conditional intensities $\lambda_n(t)$ and $\mu_n(t)$. We assume that $\lambda_n \in \mathcal{L}_n = \{\lambda_{n,\theta} : \theta \in \Theta\}$ and $\mu_n \in \mathcal{M}_n = \{ \mu_{n,\tau} : \tau \in \mathcal{T} \}$, where both parameter spaces $\Theta$ and $\mathcal{T}$ are open subsets of $\mathbb{R}^m$, so that both $\theta$ and $\tau$ are $m$-dimensional parameter vectors. 
It follows that there are $\theta_0 \in \Theta$ and $\tau_0 \in \mathcal{T}$ such that $\lambda_n = \lambda_{n ,\theta_0}$ and $\mu_n = \mu_{n,\tau_0}$. 
We assume that both conditional intensities, and the parametric families they belong to, satisfy Assumptions A1--A3 of the previous section. 
With $\widehat\theta$ and $\widehat\tau$ denoting the MLEs of $\theta_0$ and $\tau_0$, we consider the compensated processes 
\begin{equation}\label{W.hat.emp}
\begin{split}
\widehat{W}_{\lambda_n}(t) &
= \frac{1}{\sqrt{n}}\left[N_{\lambda_n}(t)-\int_{0}^{t}\lambda_{n,\widehat{\theta}}(s)\,\mathrm{d}s\right],\ \,\\
\widehat{W}_{\mu_n}(t) &
= \frac{1}{\sqrt{n}}\left[N_{\mu_n}(t)-\int_{0}^{t}\mu_{n,\widehat{\tau}}(s)\,\mathrm{d}s\right].
\end{split}
\end{equation}
By Theorem \ref{thm:1}, these two processes converge weakly in $D([0,T])$ to the respective limits
\begin{equation}\label{W.hats}
\begin{split}
V_{B_\lambda}(t) &= W_{B_\lambda}(t)
-\int_{0}^{t}\alpha_\lambda^{\top}(s)\beta_\lambda(s)\,\mathrm{d}s\,\Gamma_\lambda^{-1}\int_{0}^{T}\alpha_\lambda(s')\,\mathrm{d}W_{B_\lambda}(s'),\\
V_{B_\mu}(t) &= W_{B_\mu}(t)
-\int_{0}^{t}\alpha^{\top}_\mu(s) \beta_\mu(s)\,\mathrm{d}s\,\Gamma_\mu^{-1}\int_{0}^{T}\alpha_\mu(s')\,\mathrm{d}W_{B_\mu}(s'),
\end{split}
\end{equation}
where $\alpha_\lambda, \beta_\lambda, \Gamma_\lambda$ are as defined in Assumptions A2--A3 above, with the subscript $\theta_0$ omitted for notational brevity. 
The objects $\alpha_\mu, \beta_\mu, \Gamma_\mu$ are defined analogously, as limits arising from $\mu_{n,\tau_0}$. We also have $B_\lambda(t) = \int_0^t \beta_\lambda(s)\,\dd s$ and $B_\mu(t) = \int_0^t \beta_\mu(s)\,\dd s$, for $t \in [0,T]$.

We can rewrite the processes in \eqref{W.hats} in a slightly more concise form, by introducing the vector-valued functions
\begin{equation}\label{q.def}
\begin{split}
q_\lambda(t) = (q_{\lambda,1}(t), \ldots, q_{\lambda,m}(t))^\top := \Gamma_\lambda^{-1/2} \alpha_\lambda(t), \qquad t \in [0,T],\\
q_\mu(t) = (q_{\mu,1}(t), \ldots, q_{\mu,m}(t))^\top := \Gamma_\mu^{-1/2} \alpha_\mu(t), \qquad t \in [0,T],
\end{split}
\end{equation}
so that \eqref{W.hats} becomes
\begin{equation}\label{W.hats2}
\begin{split}
V_{B_\lambda}(t) &= W_{B_\lambda}(t)
-\int_{0}^{t}q_\lambda^{\top}(s)\beta_\lambda(s)\,\dd s\int_{0}^{T}q_\lambda(s')\,\mathrm{d}W_{B_\lambda}(s'),\\
V_{B_\mu}(t) &= W_{B_\mu}(t)
-\int_{0}^{t}q_\mu^\top(s) \beta_\mu(s)\,\dd s \int_{0}^{T}q_\mu(s')\,\mathrm{d}W_{B_\mu}(s').
\end{split}
\end{equation}

The aim of this section is to describe a transformation from $V_{B_\lambda}(t)$ to $V_{B_\mu}(t)$, where $B_\lambda$ and $B_\mu$ are absolutely continuous with respect to each other, that is, the Radon-Nikod\'{y}m derivative $\frac{\dd B_\mu}{\dd B_\lambda} = \frac{\beta_\mu}{\beta_\lambda}$ is finite and strictly positive on $[0,T]$. In our testing procedure, the process $V_{B_\mu}(t)$ will be the ``standard'' target process that any other process $V_{B_\lambda}(t)$ will be transformed into. This will unify all parametric testing problems with an $m$-dimensional parameter into a single parametric testing problem. For now, we do not specify the exact structure of the target process $V_{B_\mu}(t)$, but we will come back to this issue at the end of the section.

In what follows, it will be convenient to work with function-indexed, rather than time-indexed, processes (\citealp{vdvaart:1996}, Ch.~2.1). 
So we consider the function space $L^2(B_\lambda)$, consisting of all functions $\varphi: [0,T] \to \mathbb{R}$ satisfying
\[
\int_0^T \varphi^2(t)\,\dd B_\lambda(t) = \int_0^T \varphi^2(t) \beta_\lambda(t) \,\dd t < \infty,
\]
and the analogously defined space $L^2(B_\mu)$. Note that $L^2(B_\lambda)$ is a Hilbert space with inner product and induced norm
\begin{align*}
\langle \varphi_{1},\varphi_{2}\rangle_{B_\lambda} 
&:= \int_{0}^{T}  \varphi_{1}(t)\varphi_{2}(t)\beta_\lambda(t)\,\mathrm{d}t,  \qquad \varphi_{1},\varphi_{2} \in L^{2}(B_\lambda),\\
\|\varphi\|_{B_\lambda} &:= \langle \varphi,\varphi\rangle^{1/2}_{B_\lambda}  
=\left(\int_0^T  \varphi^2(t)\beta_\lambda(t)\,\dd t\right)^{1/2}, \qquad \varphi \in L^{2}(B_\lambda),
\end{align*}
and $L^2(B_\mu)$ is also a Hilbert space with analogously defined inner product $\langle \cdot,\cdot\rangle_{B_\mu}$ and induced norm $\|\cdot\|_{B_\mu}$. 
We can now introduce the function-indexed version of $V_{B_\lambda}$ in \eqref{W.hats2} as follows:
\begin{equation}\label{W.hatf}
\begin{split}
V_{B_\lambda}(\varphi) &:= \int_0^T \varphi(s)\, \dd V_{B_\lambda}(s)\\
&= \int_0^T \varphi(s)\, \dd W_{B_\lambda}(s) - \int_0^T \varphi(s) q^\top_\lambda(s)\beta_\lambda(s) \,\dd s\int_0^T q_\lambda(s')\, \dd W_{B_\lambda}(s'),
\end{split}
\end{equation}
for $\varphi \in L^2(B_\lambda)$. The function-indexed process $V_{B_\mu}(\psi)$ is defined analogously for $\psi \in L^2(B_\mu)$. 
Note that we can also use the inner product notation to rewrite \eqref{W.hatf} succinctly as
\begin{equation}\label{f.i}
V_{B_\lambda}(\varphi) = W_{B_\lambda}(\varphi) - \langle \varphi, q_\lambda\rangle_{B_\lambda}^\top W_{B_\lambda}(q_\lambda),
\end{equation}
where $\langle \varphi, q_\lambda\rangle_{B_\lambda}$ denotes the vector $(\langle \varphi, q_{\lambda,1}\rangle_{B_\lambda}, \ldots, \langle \varphi, q_{\lambda,m}\rangle_{B_\lambda})^\top$, and $W_{B_\lambda}(q_\lambda)$ denotes the vector $(W_{B_\lambda}(q_{\lambda,1}), \ldots, W_{B_\lambda}(q_{\lambda,m}))^\top$. 

We also note here that if we let $\langle q_\lambda, q_\lambda\rangle_{B_\lambda}$ denote the matrix of inner products $\langle q_{\lambda,i}, q_{\lambda,j}\rangle_{B_\lambda}$ for $i,j \in \{1,\ldots,m\}$, then we have
\begin{equation}\label{q.ortho}
\langle q_\lambda, q_\lambda\rangle_{B_\lambda} 
= \int_0^T q_\lambda(t)q^\top_\lambda(t)\beta_\lambda(t)\,\dd t = \Gamma_{\lambda}^{-1/2} \int_0^T\alpha_\lambda(t)\alpha^\top_\lambda(t)\beta_\lambda(t)\,\mathrm{d}t \, \Gamma_\lambda^{-1/2}=I,
\end{equation}
with $I$ denoting the $m \times m$ identity matrix. In other words, the set of functions $\{q_{\lambda,i}: 1 \le i \le m\}$ is \emph{orthonormal} in $L^2(B_\lambda)$, so $q_\lambda$ is simply the orthonormalized version of $\alpha_\lambda$ in $L^2(B_\lambda)$. Analogously, $q_\mu$ is the orthonormalized version of $\alpha_\mu$. This will play an important role in the exposition below.

\subsection{Transformation}
Before we tackle the problem of transforming $V_{B_\lambda}$ into $V_{B_\mu}$, we first consider the easier problem of transforming $W_{B_\lambda}$ into $W_{B_\mu}$. Let us denote
\begin{equation}\label{ell.def}
\ell(t) = \bigg( \frac{\dd B_\mu}{\dd B_\lambda}(t)\bigg)^{1/2} = \bigg( \frac{\beta_\mu(t)}{\beta_\lambda(t)}\bigg)^{1/2}, \qquad t \in [0,T].
\end{equation}
Note that for any $\psi \in L^2(B_\mu)$, we have $\ell \psi \in L^2(B_\lambda)$, with
\begin{equation*}
\langle\ell \psi_1, \ell \psi_2 \rangle_{B_\lambda} = \int_0^T \ell^2(t) \psi_1(t)\psi_2(t) \beta_\lambda(t)\,\dd t = \int_0^T \psi_1(t)\psi_2(t) \beta_\mu(t)\,\dd t = \langle \psi_1, \psi_2 \rangle_{B_\mu}.
\end{equation*}
Similarly, for any $\varphi \in L^2(B_\lambda)$, we have $\varphi/\ell \in L^2(B_\mu)$, with $\langle \varphi_1/\ell, \varphi_2/\ell \rangle_{B_\mu} = \langle \varphi_1, \varphi_2 \rangle_{B_\lambda}$. In other words, the linear operator $\psi \mapsto \ell \psi$ is a \emph{unitary} operator from $L^2(B_\mu)$ onto $L^2(B_\lambda)$, a bijection that preserves the inner product. 
This leads to the following result.

\begin{lem}\label{lem.elliso}
We have $\big\{W_{B_\lambda}(\ell \psi)\big\}_{\psi \in L^2(B_\mu)} \stackrel{d}{=} \big\{W_{B_\mu}(\psi)\big\}_{\psi \in L^2(B_\mu)}$.
\end{lem}

\begin{rem}
In Lemma~\ref{lem.elliso} and similar results below, ``$\stackrel{d}{=}$'' denotes equality in finite dimensional distributions. 
In particular, the statement of Lemma~\ref{lem.elliso} is equivalent to
\[
\big(W_{B_\lambda}(\ell \psi_1), \ldots, W_{B_\lambda}(\ell \psi_k)\big)^\top \stackrel{d}{=} \big( W_{B_\mu}(\psi_1), \ldots, W_{B_\mu}(\psi_k)\big)^\top
\]
for any $k \ge 1$ and any $\psi_1, \ldots, \psi_k \in L^2(B_\mu)$.
\end{rem}

Lemma \ref{lem.elliso} gives us a way of transforming $W_{B_\lambda}$ into $W_{B_\mu}$, by simply multiplying the argument by $\ell(\cdot)$. 
In order to accomplish the more difficult task of transforming $V_{B_\lambda}$ into $V_{B_\mu}$, we introduce the following family of linear operators mapping $L^2(B_\lambda)$ onto itself: given functions $a,b \in L^2(B_\lambda)$ with $\|a\|_{B_\lambda} = \|b\|_{B_\lambda} = 1$, we define the operator $R_{a,b}: L^2(B_\lambda) \to L^2(B_\lambda)$ as
\begin{equation}\label{eq.K}
R_{a,b}\, \varphi = \varphi-2\frac{\langle a-b,\varphi\rangle_{B_\lambda}}{\|a-b\|^2_{B_\lambda}}(a-b),\qquad
\varphi\in L^{2}(B_\lambda).
\end{equation}
The following lemma lists some important and convenient properties of the family of operators of the form $R_{a,b}$.
\begin{lem}\label{lem.K}
For any $a,b \in L^2(B_\lambda)$ with $\|a\|_{B_\lambda} = \|b\|_{B_\lambda} = 1$, the linear operator $R_{a,b}$ defined in \eqref{eq.K} satisfies the following properties:
\begin{itemize}
\item[(i)] $R_{a,b}$ is unitary, i.e., it is bijective and preserves the inner product.
\item[(ii)] $R_{a,b}$ is an involution, i.e., $R_{a,b}\circ R_{a,b}\, \varphi = \varphi$ for any $\varphi \in L^2(B_\lambda)$.
\item[(iii)] $R_{a,b}$ is self-adjoint, i.e., $\langle R_{a,b}\,\varphi_{1},\varphi_{2}\rangle_{B_\lambda}= \langle \varphi_{1},R_{a,b}\,\varphi_{2}\rangle_{B_\lambda}$ for any $\varphi_{1},\varphi_{2}\in L^{2}(B_\lambda)$.
\item[(iv)] $R_{a,b}\, a=b$ and $R_{a,b}\, b=a$.
\end{itemize}
\end{lem}

Lemma \ref{lem.K} suggests that the operator $R_{a,b}$ may be viewed as a ``reflection'' within $L^2(B_\lambda)$ that swaps $a$ and $b$, and leaves any function orthogonal to both of them unchanged.

We are now equipped to state our transformation result. 
For easier exposition, we will first consider the case $m=1$. 
That is, we will consider the (function-indexed) process $V_{B_\lambda}$ in \eqref{W.hatf}, and the analogously defined $V_{B_\mu}$, with functions $q_\lambda$ and $q_\mu$ taking values in $\mathbb{R}^m = \mathbb{R}$. 
Note that these processes will arise as the limits of the compensated point processes in \eqref{W.hat.emp}, when the underlying parameter spaces $\Theta$ and $\mathcal{T}$ are 1-dimensional. The transformation result involves the unitary operator $R_{a,b}$ with $a=q_\lambda$ and $b=\ell q_\mu$, i.e., $R_{q_\lambda,\ell q_\mu}$.
The result may be compared to Lemma~\ref{lem.elliso}.

\begin{thm}\label{thm.main.1}
Let $m=1$. We have
\begin{equation}\label{eq.main.1}
\Big\{V_{B_\lambda}(R_{q_\lambda,\ell q_\mu}\ell\psi)\Big\}_{\psi \in L^2(B_\mu)} \stackrel{d}{=} \Big\{V_{B_\mu}(\psi)\Big\}_{\psi \in L^2(B_\mu)}.
\end{equation}
\end{thm}

Theorem~\ref{thm.main.1} gives us a way of transforming the process $V_{B_\lambda}$ into $V_{B_\mu}$ in the case $m=1$. 
For general $m\ge 1$, we recursively define the unitary operators
\[
U^{(1)} = R_{q_{\lambda,1},\ell q_{\mu,1}} \qquad \text{ and } \qquad U^{(k)} = R_{q_{\lambda,k}, U^{(k-1)}\ell q_{\mu,k}} \circ U^{(k-1)}, \quad 2 \le k \le m,
\]
and state the following generalization of Theorem \ref{thm.main.1}, which is the main result of this section.
\begin{thm}\label{thm.main.2}
Let $m\geq 1$. We have
\begin{equation*}
\Big\{V_{B_\lambda}(U^{(m)}\ell\psi)\Big\}_{\psi \in L^2(B_\mu)} \stackrel{d}{=}
\Big\{V_{B_\mu}(\psi)\Big\}_{\psi \in L^2(B_\mu)}.
\end{equation*}
\end{thm}

For applications, it will be useful to have an analog of Theorem \ref{thm.main.2} for \emph{time-indexed} stochastic processes $\{V_{B_\lambda}(t)\}_{t \in [0,T]}$ and $\{V_{B_\mu}(t)\}_{t \in [0,T]}$. But this can be easily achieved by applying Theorem \ref{thm.main.2} to the collection of functions $\psi_t(\cdot) = \mathds{1}_{[0,t]}(\cdot), \, 0 \le t \le T,$ which leads to the following corollary.

\begin{cor}\label{cor.main.2}
For $m\geq 1$, we have
\begin{equation*}
\bigg\{\int_{0}^{T}U^{(m)}\ell\mathds{1}_{[0,t]}(s)\,\mathrm{d}V_{B_\lambda}(s) \bigg\}_{t \in [0,T]}\stackrel{d}{=}
\Big\{V_{B_\mu}(t)\Big\}_{t \in [0,T]}.
\end{equation*}
\end{cor}

\subsection{Target process}\label{sec.target}
Corollary \ref{cor.main.2} describes a transformation from $V_{B_\lambda}$ to $V_{B_\mu}$, both limiting processes as described in \eqref{W.hats2}, with the same dimension of the underlying parameter spaces. For applications, we need to specify the target process $V_{B_\mu}$, that is, explicitly describe the functions $\beta_\mu$ and $q_\mu$. We now turn to the question of how to choose this target process.

As mentioned at the end of Section \ref{sec:point}, we would like to take $V_{B_\mu}$ as the limiting form of 
\begin{equation}\label{W.mu.hat}
\widehat{W}_{\mu_n}(t) = \frac{1}{\sqrt{n}}\left[N_{\mu_n}(t)-\int_{0}^{t}\mu_{n,\widehat{\tau}}(s)\,\mathrm{d}s\right],
\end{equation}
where $N_{\mu_n}(t)$ is a Poisson process on $[0,T]$ with constant intensity $\mu_n$, and the parametric model $\mathcal{M}_n = \{\mu_{n,\tau}(t) : \tau \in \mathcal{T}\}$ includes this Poisson process for a particular value of $\tau_0 \in \mathcal{T}$. So we propose to take $\mu_n(t) = n/T$, and to embed $\mu_n$ in the parametric collection of conditional intensities
\begin{equation}\label{eq:Poissonian}
\mu_{n,\tau}(t) = \frac{n}{T} \sum_{j=1}^m \tau_j p_j(t/T), \qquad t \in [0,T],
\end{equation}
where $p_1, \ldots, p_m$ are functions mapping $[0,1]$ into $\mathbb{R}$, satisfying $p_1 \equiv 1$ and
\[
\int_{0}^{1}p_{j}(s)\,p_{k}(s)\,\mathrm{d}s=
\begin{cases} 0 &\text{ if } j \neq k \\ 1 &\text{ if } j=k\end{cases}, \qquad j,k \in \{1,\ldots,m\}.
\]
In other words, $p_1, \ldots, p_m$ are \emph{orthonormal} functions on $[0,1]$ with $p_1 \equiv 1$. 
Below we will see why it is advantageous to require $p_1, \ldots, p_m$ to be orthonormal. It is clear that $\mu_n \in \mathcal{M}_n$ with parameter $\tau_0 = (1,0,\ldots,0)^\top \in \mathbb{R}^m$. 
The parameter space $\mathcal{T}$ can be taken as any open neighborhood of $\tau_0$ small enough to ensure $\mu_{n,\tau} > 0$ on $[0,T]$, for all $\tau \in \mathcal{T}$.

Various well-known families of orthonormal functions on $[0,1]$ can be substituted for $p_1,\ldots,p_m$. 
In our applications, we will use normalized versions of the orthogonal Legendre polynomials
\vspace{-0.1cm}
\[
\widetilde{p}_k(s) = (-1)^{k-1} \sum_{j=0}^{k-1} \binom{k-1}{j} \binom{k+j-1}{j}(-s)^j, \qquad s \in [0,1], \quad k=1,\ldots,m,
\]
\vspace{-0.1cm}
that is, we will use
\vspace{-0.1cm}
\[
p_k(s) = \frac{\widetilde{p}_k(s)}{\int_0^1 \widetilde{p}_k(s')\,\dd s'}, \qquad s \in [0,1], \quad k=1,\ldots,m.
\]
\vspace{-0.1cm}
The explicit formulas for $p_1, p_2, p_3$ are:
\vspace{-0.1cm}
\[
p_1(s) = 1, \qquad p_{2}(s)=\sqrt{12}\Big(s-\frac{1}{2}\Big),\qquad
p_{3}(s)=6\sqrt{5}\Big[\Big(s-\frac{1}{2}\Big)^{2}-\frac{1}{12}\Big].
\]

We can now explicitly describe the target process we will use in our applications: it will be the process $V_{B_\mu}$ appearing in \eqref{W.hats}, with
\begin{align*}
\beta_\mu(t) &= \lim_{n \to \infty} \frac{\mu_{n,\tau_0}(t)}{n} = \frac{1}{T}, \qquad B_\mu(t) = \int_0^t \beta_\mu(s)\,\dd s = \frac{t}{T},\\
\alpha_\mu(t) &= \lim_{n \to \infty} \frac{\dt{\mu}_{n,\tau_0}(t)}{\mu_{n,\tau_0}(t)} = (p_1(t/T), \ldots, p_m(t/T))^\top =: p(t/T),\\
\Gamma_\mu &= \int_0^T \alpha^\top_\mu(t) \alpha_\mu(t)\beta_\mu(t)\,\dd t = \frac{1}{T} \int_0^T p^\top(t/T) p(t/T) \, \dd t = I.
\end{align*}
We see here the advantage of requiring $p_1, \ldots, p_m$ to be orthonormal: it gives us $\Gamma_\mu=I$, so that the function $q_\mu(\cdot) = R^{-1/2}\alpha_\mu(\cdot)$ is simply identical to $\alpha_\mu(\cdot) = p(\cdot/T)$. So we can also describe our target process as $V_{B_\mu}$ in \eqref{W.hats2} (or its function-indexed form \eqref{W.hatf}), with
\begin{equation}\label{q.mu}
\beta_\mu(t) = \frac{1}{T}, \qquad B_\mu(t) = \frac{t}{T}, \qquad q_\mu(t) = p(t/T),\qquad \mathrm{for}\ t \in [0,T]. 
\end{equation}

For any given $m \ge 1$, it is straight-forward to generate paths from the target process $V_{B_\mu}$ described above, and therefore to infer distributions of specific functionals of this process. Our next task is to describe an empirical version of the transformation in Theorem \ref{thm.main.2}, which can be applied directly to the observed path $\widehat{W}_n$ in \eqref{eq:Wnparametric}, and to show that the resulting empirical process converges to $V_{B_\mu}$ under the null hypothesis. Then functionals of this transformed empirical process will converge to the same functionals computed from $V_{B_\mu}$, which will form the basis of our goodness-of-fit testing approach. We formalize this approach in the next section.

\section{Goodness-of-Fit Testing}\label{sec:gof} 

Consider a situation where we have observed a single path of a point process $N_n(t), t \in [0,T]$, with unknown conditional intensity $\lambda_n(t)$. We have a parametric family of conditional intensities $\mathcal{L}_n = \{\lambda_{n,\theta} : \theta \in \Theta\}$, and we would like to test the null hypothesis $\lambda_n \in \mathcal{L}_n$ against the alternative $\lambda_n \notin \mathcal{L}_n$. 

Based on the results in the preceding sections, a natural approach would be the following. 
Assume that the null hypothesis is true, so that $\lambda_n = \lambda_{n,\theta_0}$ for some $\theta_0 \in \Theta$. 
Compute the MLE $\widehat{\theta}$ of $\theta_0$ from the observed path, and construct the process $\widehat{W}_n$ in \eqref{eq:Wnparametric}. 
By Theorem~\ref{thm:1}, assuming $n$ is large enough, the observed path of $\widehat{W}_n$ should behave like an observation from the limiting process $V_{B_\lambda}$. Therefore, if we apply the transformation of Corollary \ref{cor.main.2} to $\widehat{W}_n$, we should obtain a process 
\begin{equation}\label{eq.tr}
\mathscr{T}\widehat{W}_n(t) := \int_{0}^{T}U^{(m)}\ell\mathds{1}_{[0,t]}(s)\,\mathrm{d}\widehat{W}_n(s), \qquad t \in [0,T],
\end{equation}
that behaves like the target process $V_{B_\mu}$ under the assumption of the null hypothesis. 
So we can test the null hypothesis by comparing the path $\mathscr{T}\widehat{W}_n$ with the statistical behavior of $V_{B_\mu}$ using appropriate functionals.

The approach described above has the problem that the transformation $\mathscr{T}$ in \eqref{eq.tr} depends on the unknown $\theta_0$, since the function $\ell$ and the unitary operator $U^{(m)}$ both depend on $\theta_0$, and therefore the process $\mathscr{T}\widehat{W}_n$ in \eqref{eq.tr} is not constructable from the observed path of $N_n$. 
But we can remedy this problem by replacing $\theta_0$ by its estimator $\widehat{\theta}$ in $\mathscr{T}$. 
To this end, we introduce the function
\[
\widehat\ell(t) = \bigg( \frac{\beta_\mu(t)}{\widehat\beta_{\lambda}(t)} \bigg)^{1/2}, \qquad t \in [0,T],
\]
with $\widehat\beta_\lambda := \beta_{\lambda, \widehat\theta}$, as well as the unitary operators
\[
\widehat{U}^{(1)} = \widehat{R}_{\widehat{q}_{\lambda,1},\widehat{\ell} q_{\mu,1}} \qquad \text{ and } \qquad
\widehat{U}^{(k)} = \widehat{R}_{\widehat{q}_{\lambda,k}, \widehat{U}^{(k-1)} \widehat{\ell} q_{\mu,k}} \circ \widehat{U}^{(k-1)}, \quad 2 \le k \le m,
\]
where $\widehat{q}_\lambda$ and $\widehat{R}_{a,b}$ are as defined in \eqref{q.def} and \eqref{eq.K}, respectively, with the true parameter value $\theta_0$ replaced by the estimator $\widehat\theta$. 
We can then consider the following empirical analog of the process $\mathscr{T}\widehat{W}_n$ in \eqref{eq.tr}:
\begin{equation}\label{eq.trhat}
\widehat{\mathscr{T}}\,\widehat{W}_n(t) = \int_{0}^{T} \widehat{U}^{(m)} \widehat{\ell}\mathds{1}_{[0,t]}(s) \,\mathrm{d}\widehat{W}_n(s), \qquad t \in [0,T].
\end{equation}
Unlike the process $\mathscr{T}\widehat{W}_n$ in \eqref{eq.tr}, the process $\widehat{\mathscr{T}}\,\widehat{W}_n$ in \eqref{eq.trhat} \emph{can} be constructed from the observed path of $N_n$, and thus can be compared with the target process $V_{B_\mu}$. It remains to formally establish that the process $\widehat{\mathscr{T}}\,\widehat{W}_n$ indeed converges to the target process $V_{B_\mu}$. We do this under some additional technical assumptions.

For notational brevity, let us denote
\[
\Delta^U_n(s,t) := \widehat{U}^{(m)}\widehat{\ell}\mathds{1}_{[0,t]}(s) - U^{(m)}\ell\mathds{1}_{[0,t]}(s), \quad (s,t) \in [0,T]^2,
\]
and for any function $\varphi: [0,T] \to \mathbb{R}$, let $\mathrm{V}_0^T(\varphi)$ denote the \emph{total variation} of $\varphi$ over the interval $[0,T]$, i.e., the supremum of $\sum_{i=1}^n |\varphi(x_i)-\varphi(x_{i-1})|$ over all partitions $0=x_0 < x_1 < \ldots < x_n = T$. We state the following set of assumptions:

\bigskip

\textbf{A4.} For any $S \in (0,T)$, and any $m \ge 1$, we have
\vspace{10pt}
\begin{enumerate}
\begin{minipage}{0.45\linewidth}
\item[(i)] $\displaystyle \sup_{t \in [0,S]}|U^{(m)}\ell\mathds{1}_{[0,t]}(T)| < \infty$,
\item[(iii)] $\displaystyle \sup_{t \in [0,S]} |\Delta^U_n(T,t)| \stackrel{P}{\to} 0$,
\end{minipage}
\begin{minipage}{0.45\linewidth}
\item[(ii)] $\displaystyle \sup_{t \in [0,S]}\mathrm{V}_0^T(U^{(m)}\ell\mathds{1}_{[0,t]}) < \infty$,
\item[(iv)] $\displaystyle \sup_{t \in [0,S]} \mathrm{V}_0^T(\Delta^U_n(\cdot,t)) \stackrel{P}{\to} 0$.
\end{minipage}
\end{enumerate}

\bigskip

\begin{rem}
The conditions of A4 will be satisfied under appropriate differentiability assumptions on the underlying functions $\beta_{\lambda,\theta}, q_{\lambda,\theta}$ and $\beta_\mu, q_\mu$. We discuss the case $m=1$ in Appendix B.
\end{rem}

The following is the main result of this section and formalizes our goodness-of-fit approach.
\begin{thm}\label{thm.main}
Under Assumptions A1--A4, the process $\widehat{\mathscr{T}}\,\widehat{W}_n$ defined in \eqref{eq.trhat} converges weakly in $D([0,S])$ to $V_{B_\mu}$, for any $S \in (0,T)$.
\end{thm}

\begin{rem}
The fundamental result of Theorem~\ref{thm.main} is analogous to Theorem~\ref{thm:1} in the sense that it describes the asymptotic behavior of an empirical process constructable from the data.
The asymptotic behavior of $\widehat{\mathscr{T}}\,\widehat{W}_n$ is, however, ``standard'', unlike that of $\widehat{W}_n$.
Hence, statistical tests for the null hypothesis based on $\widehat{\mathscr{T}}\,\widehat{W}_n$ will be asymptotically distribution-free.
Many different functionals of the path of $\widehat{\mathscr{T}}\,\widehat{W}_n$ can be used for this purpose.
Thus, $\widehat{\mathscr{T}}\,\widehat{W}_n$ is a suitable testing process that can serve to generate a plenitude of asymptotically distribution-free tests.
\end{rem}

\section{Monte Carlo Simulations}\label{sec:mc}

In this section we study the finite-sample performance of our goodness-of-fit approach, by implementing it on simulated paths of point processes commonly used in applications. 
In Sections~\ref{sec:Aalen}--\ref{sec.software}, we consider Aalen-type survival processes and software reliability models, respectively.
In Section~\ref{sec:power}, we perform power calculations under the alternative hypothesis. 
In Appendix C, we give some details about the maximum likelihood estimation for each model. 
In Appendix D, we present additional Monte Carlo results for Aalen-type survival processes with random censoring, and for mixture cure models, both under the null, and for software reliability models under the alternative.

\subsection{Aalen-Type Survival Processes}\label{sec:Aalen}
We consider the multiplicative conditional intensity model of \cite{aalen:1978}, widely used in survival and duration analysis; see e.g., \cite{andersen:etal:1993} and \cite{kalbfleisch:2002} and the references therein. 
In its simplest form, the Aalen conditional intensity process is given by
\begin{equation}
\lambda_{n,\theta}(t)=\varsigma_{\theta}(t)[n - N_n(t-)],\qquad 0\leq t\leq T.
\label{eq:Aalen}
\end{equation}
Here, $n$ denotes the initial size of a human or animal population, $N_n(t)$ is the number of deaths observed in this population up to (and including) time $t$, and $N_{n}(t-)$ is the left limit of $N_{n}(s)$ as $s\uparrow t$, i.e., the number of deaths up to, but not including, time $t$. 
The function $\varsigma_{\theta}(t)$ is the hazard (or failure) rate, also known as the force of mortality, associated with a parametric lifetime distribution $F_\theta$, with corresponding density $F'_\theta = f_\theta$. In other words, $\varsigma_\theta(t) = f_\theta(t)/[1-F_\theta(t)]$, $t \ge 0$.

If the $n$ individuals in the initial population have i.i.d.\ lifetimes $T_1, \ldots, T_n$ generated from the lifetime distribution $F_\theta$, and there is no immigration into or emigration out of the population, then the number of deaths up to time $t$, that is, $N_n(t) = \sum_{i=1}^n \mathds{1}_{\{T_i \le t\}}$, $\quad 0 \le t \le T$, will be a point process with conditional intensity function as specified in (\ref{eq:Aalen}). To see why, consider the individual point processes $N^i(t) = \mathds{1}_{\{T_i \le t\}}$ for $i=1,\ldots,n$, with corresponding conditional intensities
\begin{align*}
\lambda_{\theta}^{i}(t)
&=\lim_{\Delta\downarrow 0}\frac{\mathbb{E}\left[N^{i}(t+\Delta)-N^{i}(t)\,|\,\mathcal{F}_{t-}\right]}{\Delta}
=\varsigma_{\theta}(t)\mathds{1}_{\{T_{i}\geq t\}}.
\end{align*}
Since we have $N_n(t) = \sum_{i=1}^n N^i(t)$, the conditional intensity of $N_n$ is equal to
\begin{align*}
\lambda_{n,\theta}(t)
&=\lim_{\Delta\downarrow 0}\frac{\mathbb{E}\left[N_{n}(t+\Delta)-N_{n}(t)\,|\,\mathcal{F}_{t-}\right]}{\Delta}
=\sum_{i=1}^{n}\lim_{\Delta\downarrow 0}\frac{\mathbb{E}\left[N^{i}(t+\Delta)-N^{i}(t)\,|\,\mathcal{F}_{t-}\right]}{\Delta}\\
&=\varsigma_{\theta}(t)\sum_{i=1}^{n}\mathds{1}_{\{T_{i}\geq t\}}
=\varsigma_{\theta}(t)\left[n-N_{n}(t-)\right],
\end{align*}
as in (\ref{eq:Aalen}).
In practical applications, not all lifetimes $T_1, \ldots, T_n$ will be observed in full, due to some individuals moving out of the population at random times. 
This results in right-censored lifetime observations, which we can also deal with, but we defer this analysis to Appendix D.

To illustrate our good\-ness-of-fit testing approach in the context of Aalen-type survival processes, we consider the point process $N_n$ with conditional intensity $\lambda_{n,\theta}$ as given in \eqref{eq:Aalen}, with force of mortality
\begin{equation}\label{weibull.mort}
\varsigma_\theta(t) = \frac{\theta_2}{\theta_1} \bigg( \frac{t_0+t}{\theta_1}\bigg)^{\theta_2-1}, \quad t \ge 0,
\end{equation}
for some fixed $t_0 \ge 0$ and a parameter vector $\theta = (\theta_1, \theta_2)^\top \in (0,\infty)^2$. 
Note that (\ref{weibull.mort}) describes the force of mortality of a Weibull excess lifetime distribution beyond the age of $t_0$, with corresponding cdf
\begin{equation}\label{weibull.cdf}
F_\theta(t) = 1 - \exp\bigg\{\bigg(\frac{t_0}{\theta_1}\bigg)^{\theta_2} - \bigg(\frac{t_0 + t}{\theta_1}\bigg)^{\theta_2}\bigg\}, \quad t \ge 0.
\end{equation}

One can imagine a population of $n$ individuals with i.i.d.\ Weibull lifetimes, each $t_0$ years old at time $t=0$. If we let $T_i$ denote the time of death (equivalently: excess lifetime after $t_0$) for individual $i$, then the counting process $N_n(t)=\sum_{i=1}^{n}\mathds{1}_{\{T_{i}\le t\}}$, keeping track of the number of deaths observed up to (and including) time $t$, will have the conditional intensity \eqref{eq:Aalen} with $\varsigma_\theta$ as specified in (\ref{weibull.mort}). 

In order to apply our testing procedure in this setup, we generate 1000 independent sample paths from the point process $N_n(t)$ as described above, on the interval $0 \le t \le T$, with $T=50$. 
We take $n=1000$, $t_0 =50$, and $\theta_1 = 86$, $\theta_2=9$. 
The values for $\theta_1$ and $\theta_2$ are chosen so that the resulting Weibull distribution is close to the empirical lifetime distributions (measured in years) typically observed in developed countries. 

From each generated sample path $N_n$, we compute the maximum likelihood estimates $\widehat\theta_1$ and $\widehat\theta_2$; see Appendix C for details. Then we construct the process $\widehat{W}_n$ described in \eqref{eq:Wnparametric}, and transform it into the process $\widehat{\mathscr{T}}\widehat{W}_n$ described in \eqref{eq.trhat}. The functions $\alpha_\theta$ and $\beta_\theta$ of \eqref{eq:alpha-beta} in this case take the form:
\[
\alpha_\theta(t) = \bigg(\!\! -\frac{\theta_2}{\theta_1}, \, \frac{1}{\theta_2} + \log \bigg( \frac{t_0 + t}{\theta_1}\bigg)\bigg)^\top, \qquad \quad \beta_\theta(t) = f_\theta(t).
\]
To implement the empirical transformation $\widehat{\mathscr{T}}$, we use the functions $\beta_{\widehat\theta}$ and $q_{\widehat\theta} = \Gamma_{\widehat\theta}^{-1/2}\alpha_{\widehat\theta}$, as well as the functions $\beta_\mu$ and $q_\mu$ described in \eqref{q.mu}.

According to Theorem \ref{thm.main}, the transformed process $\widehat{\mathscr{T}}\,\widehat{W}_n$ has the same asymptotic distribution as $\widehat{W}_{\mu_n}$ in \eqref{W.mu.hat}, which is the distribution of the target process $V_{B_\mu}$ described in Section \ref{sec.target}. 
In order to verify this, we generate 5000 independent paths from the homogeneous Poisson process $N_{\mu_n}(t)$, $0 \le t \le 50$, with constant intensity $\mu_n(t) = n/50$, and construct the corresponding paths of $\widehat{W}_{\mu_n}$ (with MLE $\widehat{\tau}$). 
To compare the distribution of the 1000 paths of $\widehat{\mathscr{T}}\,\widehat{W}_n$ with that of the 5000 paths of $\widehat{W}_{\mu_n}$, we use the following three test statistics, analogous to the well-known Kolmogorov-Smirnov, Cram\'er-von Mises, and Anderson-Darling statistics in empirical process theory:
\begin{align}\label{3stats}
\mathrm{KS}(\varphi) := \sup_{t \in [0,T]}| \varphi(t)|, \quad
\mathrm{CvM}(\varphi) := \frac{1}{T}\int_0^T \varphi^2(t) \, \mathrm{d}t, \quad \mathrm{AD}(\varphi) := \int_0^T \frac{\varphi^2(t)}{t} \, \mathrm{d}t.
\end{align}

\begin{figure}[t]
\centering
\epsfig{file=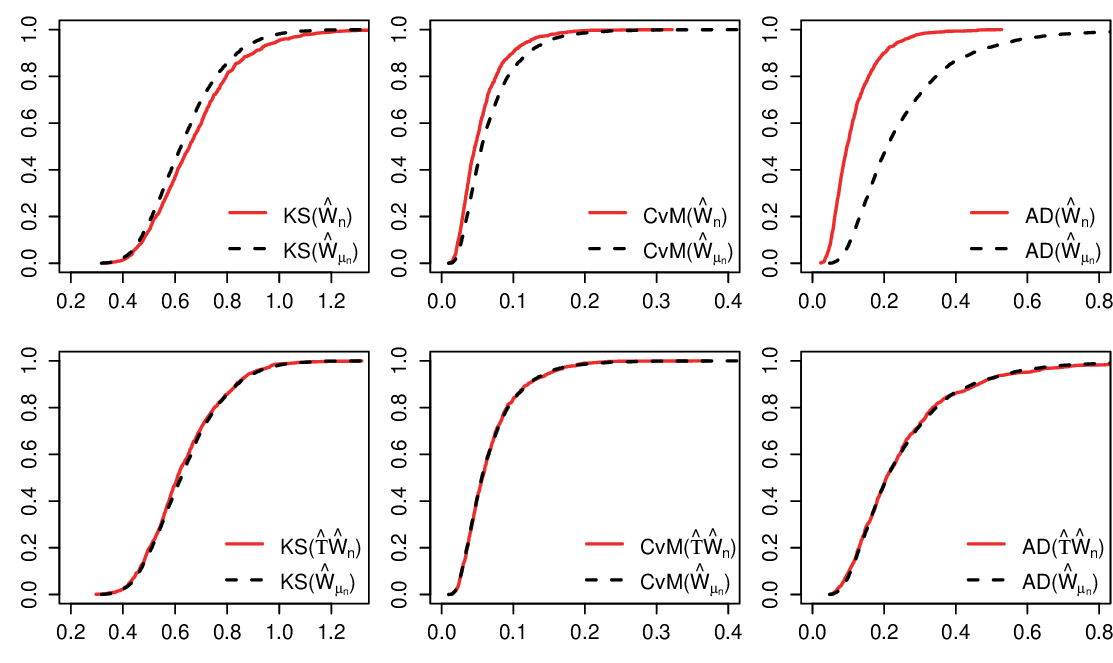,width=\textwidth}
\vspace{-15pt}
\caption{Testing for the Aalen model \eqref{eq:Aalen} under the null hypothesis. The first row shows the empirical distribution functions of the three test statistics in \eqref{3stats}, as computed from $\widehat{W}_n$ and from $\widehat{W}_{\mu_n}$. The second row shows the empirical distribution functions of the same test statistics, as computed from $\widehat{\mathscr{T}}\widehat{W}_n$ and $\widehat{W}_{\mu_n}$. Unlike the empirical process $\widehat{W}_n$, the transformed process $\widehat{\mathscr{T}}\widehat{W}_n$ seems to behave identically to the target process $\widehat{W}_{\mu_n}$, in line with Theorem \ref{thm.main}.}
\label{fig.aalen}
\end{figure}

In the first row of Fig.~\ref{fig.aalen}, we compare the empirical cdfs of these three test statistics, computed from the 1000 paths of $\widehat{W}_n$, with the empirical cdfs of the same test statistics computed from the 5000 paths of $\widehat{W}_{\mu_n}$. 
Not surprisingly, the cdfs do not match well, the process $\widehat{W}_n$ having been constructed from the Aalen-type process $N_n$ described above, and the process $\widehat{W}_{\mu_n}$ from a time-homogeneous Poisson process $N_{\mu_n}$ embedded in an ``artificial'' parametric family as described in Section~\ref{sec.target}. 
In the second row of Fig.~\ref{fig.aalen}, we compare the empirical cdfs of the three test statistics, now computed from the 1000 paths of the \emph{transformed} process $\widehat{\mathscr{T}}\,\widehat{W}_n$, with the empirical cdfs of the same test statistics computed from the 5000 paths of $\widehat{W}_{\mu_n}$. 
In line with Theorem~\ref{thm.main}, we see that the empirical distributions are now almost identical for each test statistic, confirming that the transformed process $\widehat{\mathscr{T}}\,\widehat{W}_n$ originating from the Aalen-type point process behaves just like $\widehat{W}_{\mu_n}$ originating from a homogeneous Poisson process.

Table~\ref{table.Aalen} shows the observed rejection rates at 10\%, 5\% and 1\% significance levels for the three test statistics. 
More precisely, for each test statistic we report the number of observed values that exceed the 90\%, 95\% and 99\% quantiles of the null distribution as computed from a homogeneous Poisson process. 
The observed numbers are close to the expected numbers at all significance levels, in line with the good match we observe in the second row of Fig.~\ref{fig.aalen}.

\begin{table}[t]
\begin{tabular}{c|cccc}
\textbf{Sig.\ level}   & \textbf{KS} & \textbf{CvM} & \textbf{AD} & \textbf{Expected} \\[2pt]
\hline\\[-7pt]
0.10 & 95 & 102 & 112 & 100\\[2pt]
\hline\\[-7pt]
0.05 & 44 & 55 & 54 & 50\\[2pt]
\hline\\[-7pt]
0.01 & 10 & 8 & 16 & 10\\[2pt]
\hline
\end{tabular}
\vspace{10pt}
\caption{Testing for the Aalen model \eqref{eq:Aalen} under the null hypothesis. The table shows the number of times each test statistic lead to a rejection of the null hypothesis, among 1000 simulated paths. The observed rejection rates are close to the significance levels, in line with Theorem \ref{thm.main}.}\label{table.Aalen}
\vspace{-10pt}
\end{table}

\subsection{Software Reliability Models}\label{sec.software}
In software reliability theory, point processes with conditional intensities of the type
\begin{equation}\label{eq:soft}
\lambda_{n,\theta,p}(t) = \varsigma_\theta(t)[np-N_n(t-)], \quad 0 \le t \le T,
\end{equation}
are widely used to model the occurrence of software failures over time. 
Given a software of size $n$ (where ``size'' can refer to, e.g., lines of code), assume that there is a fixed, unknown number $np$ of faults or ``bugs'' initially present in the software, proportional to software size. 
Also assume that each of these faults will lead to a software failure at a random time in the future, and the times until failure are i.i.d.\ random variables with cdf $F_\theta$ and failure rate $\varsigma_\theta$. 
If $N_n(t)$ denotes the number of software failures detected up to and including time $t$, then $N_n(t)$ will be a point process with conditional intensity \eqref{eq:soft}, by analogous arguments as in Section~\ref{sec:Aalen}. 
These software reliability models are formally related to the mixture cure models analyzed in Appendix D.
If the proportion $p$ of initial faults, and therefore the initial number $np$ of faults, is known, then the model \eqref{eq:soft} becomes equivalent to the Aalen model \eqref{eq:Aalen}.

Picking $\varsigma_\theta(t) = \theta$ for some fixed $\theta>0$ leads to the so-called Jelinski-Moranda model, first proposed in \cite{jelinski:moranda:1972}, and still widely used in software reliability analysis despite its relative simplicity (see, e.g., \cite{lyu:1996} and \cite{turk:alsolami:2016}):
\begin{equation}\label{JM}
\lambda_{n,\theta,p}(t) = \theta (np - N_n(t-)), \quad 0 \le t \le T.
\end{equation}
Note that the constant failure rate $\varsigma_\theta(t) = \theta$ implies i.i.d.\ Exp($\theta$) waiting times until failure, with cdf $F_\theta(t) = 1-\exp(-\theta t), \; t \ge 0$.

A generalization of the Jelinski-Moranda model is provided by \cite{littlewood:1980}, which assumes the following conditional intensity for the occurrence of software failures:
\begin{equation}\label{Lw}
\lambda_{n,\theta,p}(t) = \frac{\theta_1}{1 + \theta_2t}(np - N_n(t-)), \quad 0 \le t \le T.
\end{equation}
In contrast to the Jelinski-Moranda model, where each software failure decreases the conditional intensity of future failures by a fixed amount $\theta$, the Littlewood model assumes that the decrease in the conditional intensity becomes smaller with time. 
The idea is that ``obvious'' or large faults will be discovered and fixed rather quickly, while more ``subtle'' faults will take a long time to discover and fix. Note that the Littlewood model (\ref{Lw}) reduces to the Jelinski-Moranda model (\ref{JM}) as $\theta_2 \downarrow 0$. 
Also note that the decreasing failure rate $\varsigma_\theta(t) = \frac{\theta_1}{1 + \theta_2 t}, \; t \ge 0$ corresponds to i.i.d.\ Lomax($\theta_1/\theta_2$, $1/\theta_2$) waiting times until failure, with cdf $F_\theta(t) = 1 - (1+\theta_2 t)^{-\theta_1/\theta_2}, \; t \ge 0$.

In order to test the good\-ness-of-fit of the parametric intensity model (\ref{JM}), we generate 1000 independent sample paths from the corresponding point process $N_n(t)$ on the interval $0 \le t \le T$. 
We take $T=1$, $n=10{,}000$, $\theta=1$, and $p=0.1$. 
So we imagine a software of 10{,}000 lines, initially containing 1000 undiscovered faults. From each generated path $N_n(t)$, we compute the maximum likelihood estimates $\widehat\theta$ and $\widehat{p}$, and construct the processes $\widehat{W}_n$ and $\widehat{\mathscr{T}}\,\widehat{W}_n$ as before. 
The functions $\alpha_\theta$ and $\beta_\theta$ defined in \eqref{eq:alpha-beta} take the form:
\[
\alpha_{\theta,p}(t) = \bigg(\frac{1}{\theta}, \frac{e^{\theta t}}{p}\bigg)^\top, \qquad \quad \beta_{\theta,p}(t) = p \theta e^{-\theta t}.
\]

\begin{figure}[t]
\centering
\epsfig{file=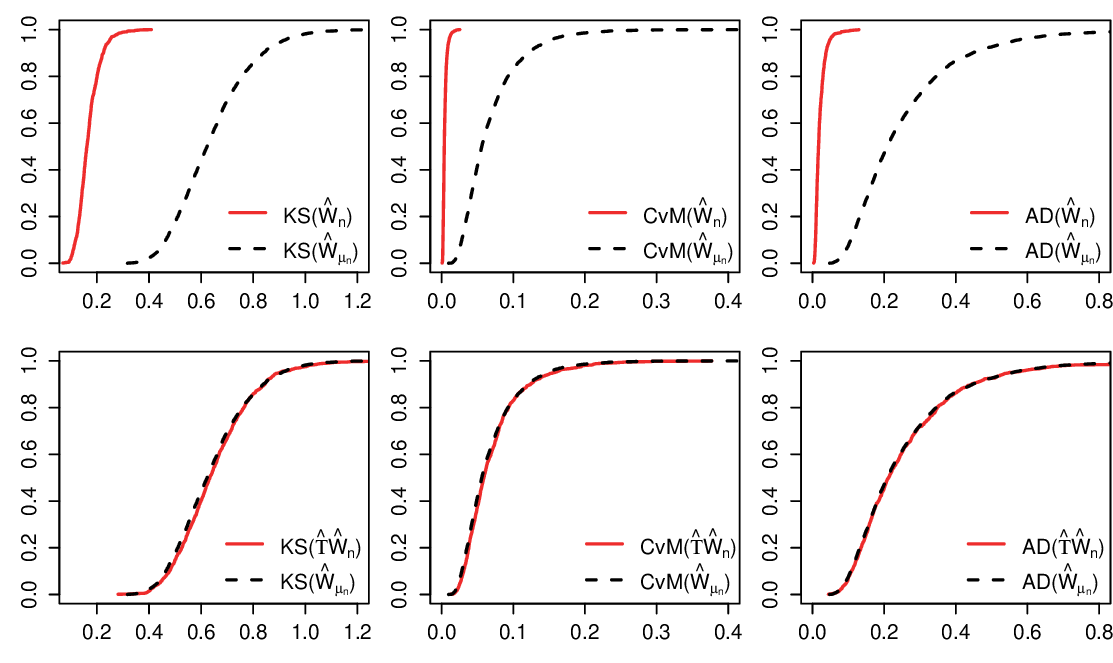,width=\textwidth}
\vspace{-15pt}
\caption{Testing for the Jelinski-Moranda model \eqref{JM} under the null hypothesis. The first row shows the empirical distribution functions of the three test statistics in \eqref{3stats}, as computed from $\widehat{W}_n$ and from $\widehat{W}_{\mu_n}$. The second row shows the empirical distribution functions of the same test statistics, as computed from $\widehat{\mathscr{T}}\widehat{W}_n$ and $\widehat{W}_{\mu_n}$. Unlike the empirical process $\widehat{W}_n$, the transformed process $\widehat{\mathscr{T}}\widehat{W}_n$ seems to behave identically to the target process $\widehat{W}_{\mu_n}$, in line with Theorem \ref{thm.main}.}
\label{fig.JM}
\end{figure}

Finally, we compare the processes $\widehat{W}_n$ and $\widehat{\mathscr{T}}\,\widehat{W}_n$ with $\widehat{W}_{\mu_n}$ using the test statistics in (\ref{3stats}). 
Fig.\ \ref{fig.JM} is the analogue of Fig.\ \ref{fig.aalen} above, and we observe once again that $\widehat{W}_n$ is pretty far from $\widehat{W}_{\mu_n}$ in statistical behavior, while $\widehat{\mathscr{T}}\,\widehat{W}_n$ is nearly identical to $\widehat{W}_{\mu_n}$, in line with Theorem~\ref{thm.main}.

Next, we also generate 1000 independent sample paths from the point process $N_n(t)$ defined by the Littlewood conditional intensity (\ref{Lw}), on the interval $[0,T]$ with $T=1$. We take $n=10{,}000$, $\theta_1=4$, $\theta_2=1$, and $p=0.1$. From each realized path, we compute the maximum likelihood estimates $\widehat\theta_1, \widehat\theta_2, \widehat p$, and construct the processes $\widehat{W}_n$ and $\widehat{\mathscr{T}}\,\widehat{W}_n$ as before. The functions $\alpha_\theta$ and $\beta_\theta$ in (\ref{eq:alpha-beta}) now take the form:
\[
\alpha_{\theta,p}(t) = \bigg(\frac{1}{\theta_1}, -\frac{t}{1+\theta_2 t}, \frac{(1 + \theta_2 t)^{\theta_1/\theta_2}}{p} \bigg)^\top, \qquad \quad \beta_{\theta,p}(t) = p\theta_1 (1+\theta_2 t)^{-1-\theta_1/\theta_2}.
\]

\begin{figure}[t]
\centering
\epsfig{file=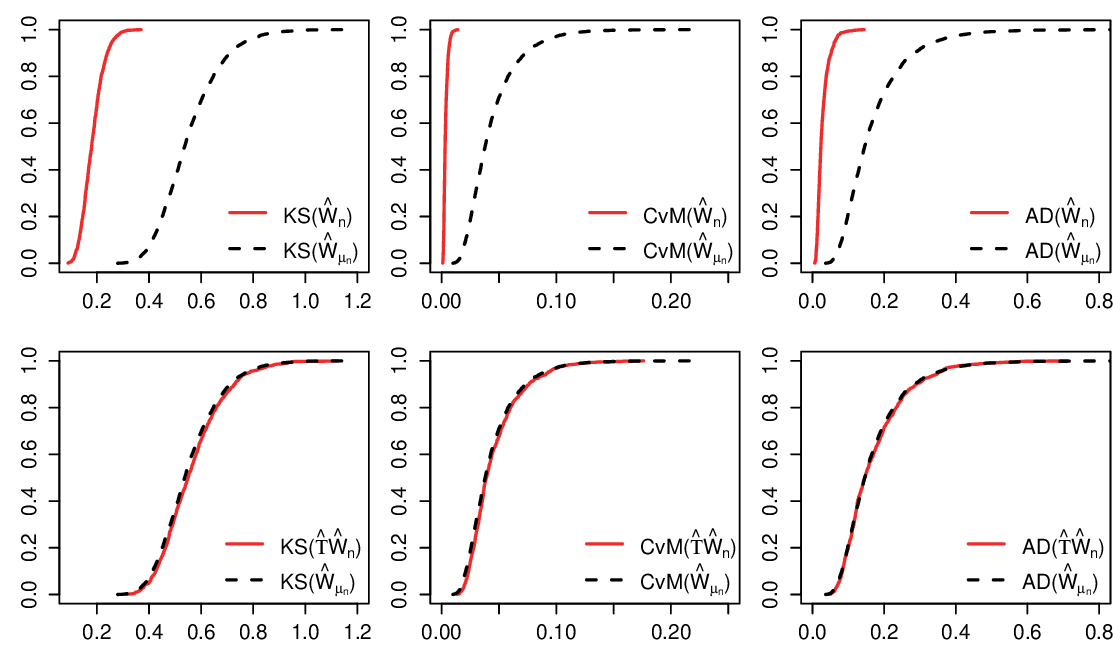,width=\textwidth}
\vspace{-15pt}
\caption{Testing for the Littlewood model \eqref{Lw} under the null hypothesis. The first row shows the empirical distribution functions of the three test statistics in \eqref{3stats}, as computed from $\widehat{W}_n$ and from $\widehat{W}_{\mu_n}$. The second row shows the empirical distribution functions of the same test statistics, as computed from $\widehat{\mathscr{T}}\widehat{W}_n$ and $\widehat{W}_{\mu_n}$. Unlike the empirical process $\widehat{W}_n$, the transformed process $\widehat{\mathscr{T}}\widehat{W}_n$ seems to behave identically to the target process $\widehat{W}_{\mu_n}$, in line with Theorem \ref{thm.main}.}
\label{fig.Lw}
\end{figure}

Then we compare the processes $\widehat{W}_n$ and $\widehat{\mathscr{T}}\,\widehat{W}_n$ with $\widehat{W}_{\mu_n}$ using the test statistics in (\ref{3stats}). Fig.\ \ref{fig.Lw} is the analogue of Fig.\ \ref{fig.aalen}--\ref{fig.JM} above, and shows similar results. 

Table~\ref{table.JMLw} is the analogue of Table \ref{table.Aalen} for the Jelinski-Moranda and Littlewood models, and shows that the observed rejection rates are once again close to their expected values.

\begin{table}[t]
\bgroup
\def\arraystretch{1.3}
\begin{tabular}{c|ccc|ccc|c}
 & \multicolumn{3}{|c|}{\textbf{Jelinski-Moranda}} & \multicolumn{3}{|c|}{\textbf{Littlewood}} & \\[2pt]
\textbf{Sig.\ level}   & \textbf{KS} & \textbf{CvM} & \textbf{AD} & \textbf{KS} & \textbf{CvM} & \textbf{AD} & \textbf{Expected} \\[2pt]
\hline
0.10 & 108 & 106 & 98 & 120 & 113 & 107 & 100\\[2pt]
\hline
0.05 & 46 & 59 & 50 & 54 & 62 & 59 & 50\\[2pt]
\hline
0.01 & 15 & 12 & 16 & 14 & 11 & 11 & 10\\
\hline
\end{tabular}
\egroup
\vspace{10pt}
\caption{Testing for the Jelinski-Moranda and Littlewood models \eqref{JM}-\eqref{Lw} under the null hypothesis. The table shows the number of times each test statistic lead to a rejection of the null hypothesis, among 1000 simulated paths per model. The observed rejection rates are close to the significance levels, in line with Theorem \ref{thm.main}.}\label{table.JMLw}
\vspace{-10pt}
\end{table}

We emphasize here that the empirical cdfs computed from $\widehat{W}_{\mu_n}$ plotted in Fig.~\ref{fig.Lw} are different from the empirical cdfs plotted in Fig.~\ref{fig.aalen}--\ref{fig.JM} above: since the process $\widehat{W}_n$ now involves \emph{three} estimated parameters rather than two, our transformation $\widehat{\mathscr{T}}$ is now a composition of three unitary transformations rather than two, and the resulting process $\widehat{\mathscr{T}}\,\widehat{W}_n$ should be compared to $\widehat{W}_{\mu_n}$ constructed with three estimated parameters rather than two. More explicitly, $\widehat{W}_{\mu_n}$ has the form \eqref{W.mu.hat} as before, but with $\mu_{n,\tau}$ in \eqref{eq:Poissonian} having $m=3$ rather than $m=2$. 

\subsection{Power Studies}\label{sec:power}

To investigate the power properties of our testing approach, we generate paths from a point process $N_n(t)$ with a conditional intensity $\lambda_{n}(t)$ that lies \emph{outside} of a particular parametric family $\mathcal{L}_n$, but construct the processes $\widehat{W}_n$ and $\widehat{\mathscr{T}}\,\widehat{W}_n$ under the incorrect assumption $\lambda_n \in \mathcal{L}_n$, and observe the behavior of the three test statistics in (\ref{3stats}). 
In this section, we do this for the Aalen-type models of Section \ref{sec:Aalen}. 
In Appendix D, we perform similar calculations for the software reliability models of Section~\ref{sec.software}, and obtain similar results.

In the case of the Aalen model, we generate 1000 independent sample paths from the process $N_n(t)$ with conditional intensity (\ref{eq:Aalen}), where $\varsigma_\theta$ is the Weibull force of mortality \eqref{weibull.mort}. 
As in Section~\ref{sec:Aalen}, we take $n=1000$, $t_0=50$, $\theta_1=86$, $\theta_2=9$. 
From each generated sample path we construct the process $\widehat{W}_n$ in \eqref{eq:Wnparametric}, now under the false assumption of Gompertz (rather than Weibull) lifetimes, with force of mortality and cdf
\begin{equation}\label{gomp.mort}
\varsigma_{\theta}(t) = \theta_{1}\theta_{2}\exp\left(\theta_{2}\left(t_{0}+t\right)\right),\quad
F_\theta(t) = 1 - \exp\{-\theta_1e^{\theta_2t_0}(e^{\theta_2 t}-1)\}, \quad t \ge 0;
\end{equation}
see Appendix C for MLE details.
\begin{figure}[t]
\centering
\epsfig{file=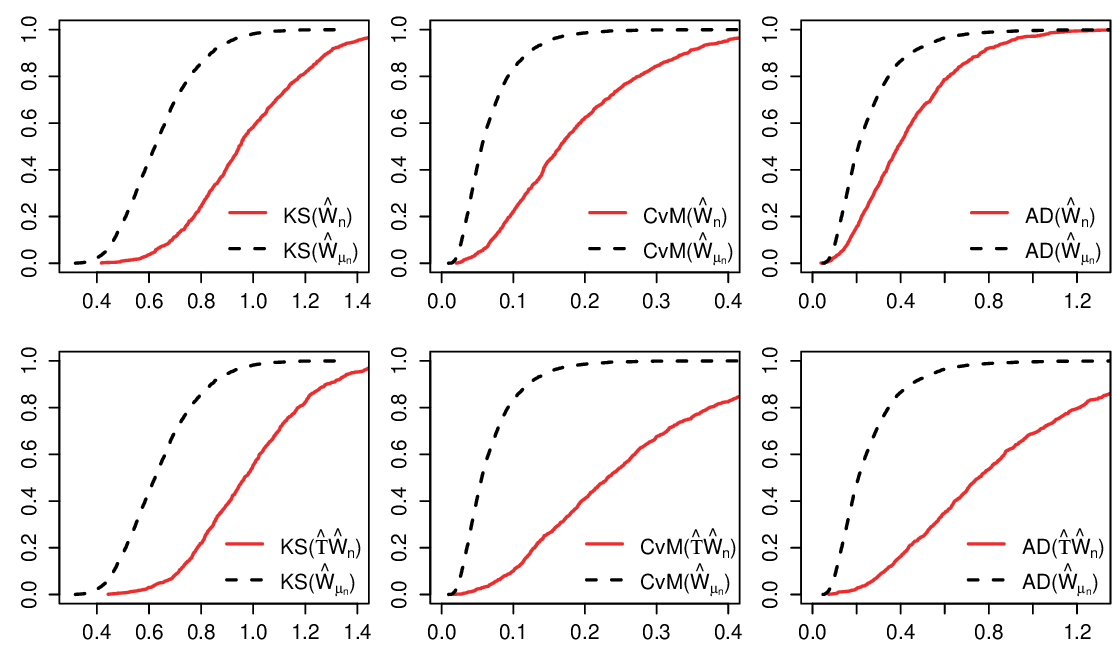,width=\textwidth}
\vspace{-15pt}
\caption{Testing for the Aalen model \eqref{eq:Aalen} under the alternative hypothesis. The first row shows the empirical distribution functions of the three test statistics in \eqref{3stats}, as computed from $\widehat{W}_n$ and from $\widehat{W}_{\mu_n}$. The second row shows the empirical distribution functions of the same test statistics, as computed from $\widehat{\mathscr{T}}\widehat{W}_n$ and $\widehat{W}_{\mu_n}$. In contrast to the situation under the null hypothesis, the transformed process $\widehat{\mathscr{T}}\widehat{W}_n$ is not any closer to the target process $\widehat{W}_{\mu_n}$ than the un-transformed $\widehat{W}_n$, implying good power properties.}
\label{fig.power1}
\end{figure}
Then we construct the transformed process $\widehat{\mathscr{T}}\,\widehat{W}_n$ in \eqref{eq.trhat}, again under the false assumption of Gompertz lifetimes. 
The functions $\alpha_{\theta}$ and $\beta_\theta$ of 
(\ref{eq:alpha-beta}) take the form
\[
\alpha_\theta(t) = \Big( \frac{1}{\theta_1}, \, \frac{1}{\theta_2} + t_0 + t \Big)^\top, \qquad \quad \beta_\theta(t) = f_\theta(t).
\]
Using the empirical cdfs of the test statistics in (\ref{3stats}) as before, we compare the un-transformed process $\widehat{W}_n$, as well as the transformed one $\widehat{\mathscr{T}}\,\widehat{W}_n$ (both computed under the misspecified Gompertz model), with $\widehat{W}_{\mu_n}$. 
See Fig.~\ref{fig.power1} for the results. 
We observe that, in contrast to the examples under the null hypothesis, in the present case the transformed process $\widehat{\mathscr{T}}\,\widehat{W}_n$ is not any closer to the target $\widehat{W}_{\mu_n}$ than the un-transformed process $\widehat{W}_n$. 
Table~\ref{table.power1} displays the observed rejection rates (out of 1000) at 10\%, 5\% and 1\% significance levels. 
As suggested by Fig.\ \ref{fig.power1}, the rejection rates are much higher than what would be expected under the null hypothesis, and confirm good power properties.

\begin{table}[b]
\begin{tabular}{c|cccc}
\textbf{Sig.\ level}   & \textbf{KS} & \textbf{CvM} & \textbf{AD} & \begin{tabular}[c]{@{}c@{}}\textbf{Expected}\\\textbf{(under null)}\end{tabular} \\[2pt]
\hline\\[-7pt]
0.10 & 713 & 834 & 789 & 100\\[2pt]
\hline\\[-7pt]
0.05 & 599 & 740 & 693 & 50\\[2pt]
\hline\\[-7pt]
0.01 & 369 & 544 & 440 & 10\\[2pt]
\hline
\end{tabular}
\vspace{10pt}
\caption{Testing for the Aalen model \eqref{eq:Aalen} under the alternative hypothesis. The table shows the number of times each test statistic lead to a rejection of the null hypothesis, among 1000 simulated paths. In contrast to the situation under the null hypothesis, the observed rejection rates are much higher than the corresponding significance levels, implying good power properties.}\label{table.power1}
\vspace{-10pt}
\end{table}

\section{Data Analysis}\label{sec:data} 

\subsection{Mortality Study}\label{sec:data:mort} We would like to apply our testing approach to a real-life data set of observed human lifetimes. 
Since lifetime data for human cohorts are rather difficult to obtain, we will instead simulate a realistic data set of lifetimes using real-life annual mortality rates reported in the Human Mortality Database, \cite{HMD}.
This is common practice in actuarial and demographic statistics.
The HMD is a collaboration between the Department of Demography at UC Berkeley, the Max Planck Institute of Demographic Research, and the French Institute for Demographic Studies, and it provides mortality data at national level for more than 40 countries, typically dating back to the late 1800s. 

For our study, we focus on the cohort of 50-year old males who were living in Luxembourg in the year 1970. 
The HMD reports that there were 1788 males of age 50 in Luxembourg at the beginning of 1970. 
We do not have the exact times of death for this cohort, but we do have \emph{death rates} recorded at ages $50, \ldots, 101$ for Luxembourgish males born in 1920 (who were 50 years old in 1970). 
The death rate $\varsigma_x$ at age $x$ is defined as the ratio of the number of deaths at age $x$ to the relevant exposure-to-risk. 
The exposure-to-risk is the estimated size of the population exposed to the risk of death, with a correction that reflects the timing of deaths during the year and the variation in the cohort's birthdays by month. 
We note here that the HMD also provides explicit numbers of deaths at ages $x=50,\ldots,101$ in years $y=1970,\ldots,2021$ recorded in Luxembourg. 
However, since these numbers include deaths among people who moved into the country after 1970, and exclude deaths among people who moved out of the country after 1970, we cannot treat them as observed numbers of death among the same initial cohort. 
Hence, we opt for simulating lifetimes based on observed mortality rates.  

\begin{figure}[b]
\centering
\epsfig{file=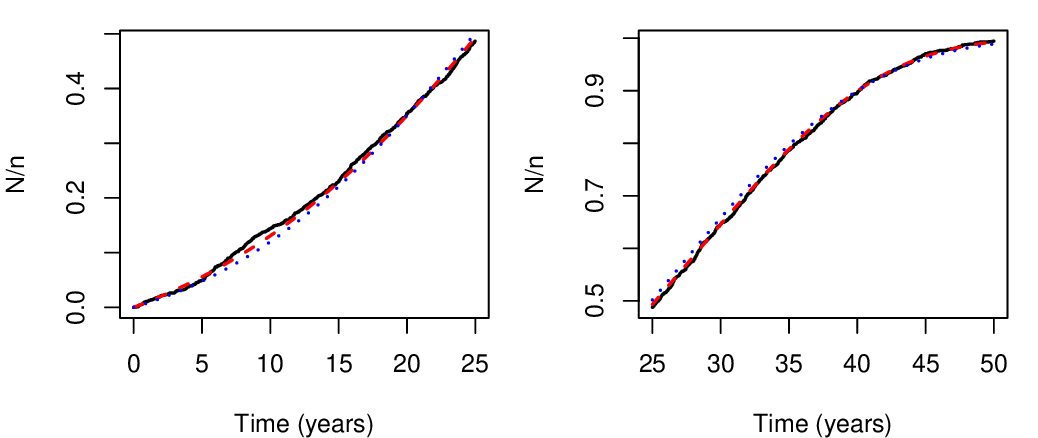,width=\textwidth}
\vspace{-8pt}
\caption{The observed process $N_n(t)/n$ (solid line) compared with its estimated limits under the Weibull model (dotted line) and the Gompertz model (dashed line). The plot is split into two periods to make the differences more visible.}
\label{fig.luxsim2}
\end{figure}

We will assume that the 1788 excess lifetimes beyond the age of 50 in our cohort are i.i.d.\ realizations from an excess lifetime distribution with piecewise constant force of mortality
\[
\varsigma(t) = \sum_{x=0}^{51} \varsigma_{50+x} \mathds{1}_{\{x \le t < x+1\}}, \quad 0 \le t < 52,
\]
with $t$ in units of years. 
This allows us to construct the cdf of the excess lifetimes as
\[
F(t) = 1 - \exp \bigg\{ \! - \int_0^t \varsigma(s)\,\dd s\bigg\} = 1 - \exp\bigg\{\! -\sum_{x=0}^{\lfloor t \rfloor -1} \varsigma_{50+x} - (t - \lfloor t \rfloor)\varsigma_{50+\lfloor t \rfloor}\bigg\},
\]
for $0 < t < 52$, and we set $F(t) = 1$ for $t \ge 52$. We can now simulate the remaining lifetimes of our cohort of 50-year old males, by randomly sampling from this cdf 1788 times. 
We thus obtain a simulated path of the point process $N_n(t)$ tracking the cumulative number of deaths over time; see Fig.~\ref{fig.luxsim2} (solid line). 
We limit our observation window to $t \in [0,T]$, with $T=50$ as in Section~\ref{sec:Aalen}. 
The total number of deaths observed in this period is 1778, which means that there are 10 people who survive beyond the age of 100 in our simulated data set. 

Our aim is to test the goodness-of-fit of the Aalen-type conditional intensity model (\ref{eq:Aalen}), with Weibull force of mortality (\ref{weibull.mort}) as well as with Gompertz force of mortality (\ref{gomp.mort}), for our simulated path $N_n$. So we follow the methodology of Section \ref{sec:Aalen} to compute maximum likelihood estimates and construct the processes $\widehat{W}_n$ and $\widehat{\mathscr{T}}\,\widehat{W}_n$ under both models. From the transformed process $\widehat{\mathscr{T}}\,\widehat{W}_n$, we compute the values of the three test statistics in (\ref{3stats}) as well as the $p$-values associated with them, using the null distributions inferred from the target process $\widehat{W}_{\mu_n}$. The results are summarized in Table~\ref{table.luxsim}. All three test statistics clearly reject the Weibull model at any reasonable significance level, while the Gompertz model is not rejected by any test statistic. We conclude that Gompertz outperforms Weibull in terms of goodness-of-fit.

\begin{table}[t]
\begin{tabular}{lcc}
\textbf{Model}   & \textbf{Parameter estimates} & \textbf{Test statistics}\\[2pt]
\hline\\[-6pt]
Weibull & $\begin{aligned} \widehat{\theta}_1 &= 78.419 \\ \widehat{\theta}_2 &= 6.242 \end{aligned}$ & $\begin{aligned} \text{KS:} & \; 1.221 \; (\mathbf{0.0010}) \\ \text{CvM:} &\; 0.302 \; (\mathbf{0.0012}) \\ \text{AD:} &\; 1.026 \; (\mathbf{0.0030}) \end{aligned}$\\[17pt]
\hline\\[-6pt]
Gompertz & $\begin{aligned} \widehat\theta_1 &= 0.003\\ 
\widehat\theta_2 &= 0.073 \end{aligned}$ & $\begin{aligned} \text{KS:} & \; 0.798 \; (0.1444) \\ \text{CvM:} &\; 0.093 \; (0.1986) \\ \text{AD:} &\; 0.373 \; (0.1642) \end{aligned}$\\[20pt]
\hline
\end{tabular}
\vspace{10pt}
\caption{Testing results for simulated lifetime data. 
The numbers in parentheses are the $p$-values corresponding to the test statistics. Bold type indicates rejection at 5\% significance.}\label{table.luxsim}
\end{table}

In Fig.~\ref{fig.luxsim2}, we compare the observed path $N_n(t)/n$ with its estimated limiting form $F_{\widehat{\theta}}(t)$ under both models; 
see (\ref{weibull.cdf}) and (\ref{gomp.mort}) for the explicit form of $F_\theta$ in the two models. While both models seem to provide a reasonably close fit to the observed path, the fit of the Weibull model is visibly worse than that of Gompertz, underestimating the numbers of deaths in years 5-20 (i.e., ages 55-70) and overestimating the numbers of deaths in years 20-40 (i.e., ages 70-90). This discrepancy is apparently picked up by our testing procedure, leading to a rejection of the Weibull model.

\subsection{Software Failures}\label{sec:data:soft} 
\begin{figure}[t]
\centering
\epsfig{file=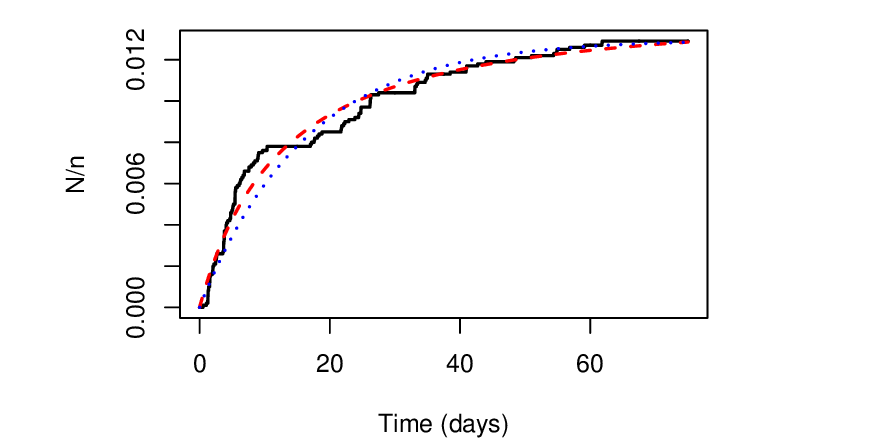,width=0.85\textwidth}
\vspace{-5pt}
\caption{The observed process $N_n(t)/n$ (solid line) compared with its estimated limits under the Jelinski-Moranda model (dotted line) and the Littlewood model (dashed line)}
\label{fig.csr2Lw}
\end{figure}

We will also apply our testing approach to a real-life data set of software failures. 
The data set comes from \cite{lyu:1996}, where it is referred to as CSR2, and it consists of 129 waiting times between successive software failures, observed over a period of 75 days. 
According to \cite{lyu:1996}, the data was sourced from a single-user workstation at the Centre for Software Reliability (CSR). 
CSR2 is part of a larger data set, CSR1, which consists of waiting times between 397 successive ``user-perceived events'' that comprise genuine software failures as well as problems related to usability and inadequate documentation. 
The cumulative number of observed failures $N_n(t)$ in CSR2 is plotted in Fig.~\ref{fig.csr2Lw} (solid line), for $0 \le t \le T$ with $T=75$ days.

We would like to test whether either of the two models discussed in Section~\ref{sec.software}, namely the Jelinski-Moranda and the Littlewood models, provides a good fit to the data set CSR2. 
To this end, we compute the maximum likelihood estimates of the parameters, construct the process $\widehat{W}_n$ and the transformed process $\widehat{\mathscr{T}}\,\widehat{W}_n$ under both models. 
Note that $n$, the ``software size'', is not known here, but this is not important: the quantity $n$ only appears as $np$ in the likelihood function and everywhere else in our procedure. As a result, scaling $n$ by some constant $c>0$ only has the effect of scaling the estimate $\widehat{p}$ by $1/c$, and leaves everything else unchanged. 
In particular, the process $\widehat{\mathscr{T}}\,\widehat{W}_n$ is independent of $n$. We arbitrarily take $n=10{,}000$ as in our numerical simulations above.

\begin{table}[b]
\begin{tabular}{lcc}
\textbf{Model}   & \textbf{Parameter estimates} & \textbf{Test statistics}\\[2pt]
\hline\\[-6pt]
Jelinski-Moranda & $\begin{aligned} \widehat\theta &= 0.061 \\ n\widehat{p} &= 129.8 \end{aligned}$ & $\begin{aligned} \text{KS:} & \; 1.199 \; (\mathbf{0.0012}) \\ \text{CvM:} &\; 0.384 \; (\mathbf{0.0004}) \\ \text{AD:} &\; 1.637 \; (\mathbf{0.0000}) \end{aligned}$\\[17pt]
\hline\\[-6pt]
Littlewood & $\begin{aligned} \widehat\theta_1 &= 0.081\\ 
\widehat\theta_2 &= 0.064 \\ n\widehat{p} &= 144.2 \end{aligned}$ & $\begin{aligned} \text{KS:} & \; 0.711 \; (0.1016) \\ \text{CvM:} &\; 0.083 \; (0.0664) \\ \text{AD:} &\; 0.665 \; (\mathbf{0.0018}) \end{aligned}$\\[20pt]
\hline
\end{tabular}
\vspace{10pt}
\caption{Testing results for CSR2. 
The numbers in parentheses are the $p$-values corresponding to the test statistics. 
Bold type indicates rejection at 5\% significance.}\label{table.csr2}
\end{table}

From the process $\widehat{\mathscr{T}}\,\widehat{W}_n$, we compute the value of the three test statistics in (\ref{3stats}), and we also compute the $p$-values associated with these observations, using the null distributions computed from the target process $\widehat{W}_{\mu_n}$. 
Note that the target process is different for the two models, due to the different numbers of parameters. 
The results are given in~Table \ref{table.csr2}. 
We see that all three test statistics reject the Jelinski-Moranda model at 5\% significance level, while the Littlewood model is not rejected by Kolmogorov-Smirnov and Cram\'er-von Mises tests, although it is rejected by Anderson-Darling. 
It appears that the Littlewood model is a more plausible model for this data set than Jelinski-Moranda.

In Fig.~\ref{fig.csr2Lw}, we compare the observed path $N_n(t)/n$ with its estimated limiting form under the Jelinski-Moranda model as well as under the Littlewood model, given by
\[
\widehat{p} F_{\widehat{\theta}}(t) = \widehat{p}\big[1-\exp(-\widehat\theta t)\big], \quad\mathrm{and}\quad \widehat{p} F_{\widehat{\theta}}(t) = \widehat{p}\big[1-(1+\widehat{\theta}_2 t)^{-\widehat\theta_1/\widehat\theta_2}\big], \quad t > 0,
\]
respectively. We observe that the Littlewood model provides a closer match to the observed path than the Jelinski-Moranda model, in agreement with our test results.

\section*{Data and Code Availability}
The mortality data used in Section \ref{sec:data:mort} is publicly available (after free registration) at: \url{https://www.mortality.org/Country/Country?cntr=LUX}. The software failure data used in Section \ref{sec:data:soft} is publicly available at: \url{https://www.cse.cuhk.edu.hk/~lyu/book/reliability/data.html}.
{\tt{R}} code to implement the procedure developed in this paper is available from the authors upon request. 

\section*{Acknowledgements}

We are very grateful to Dan Crisan, Irene Gijbels, 
David Steinsaltz and conference participants at 
the Applied Stochastic Models and Data Analysis International Conference in Florence, 
the International Symposium on Non-Parametric Statistics in Paphos, Cyprus, 
the Workshop on Statistical Methods for Dynamical Stochastic Models at Imperial College London 
and the  
International Conference on Probability Theory and Statistics at  
Tbilisi State University
for their comments and suggestions.
Roger Laeven was supported in part by the Netherlands Organization for Scientific Research under an NWO-Vici grant 2020--2027. 

\bibliographystyle{apalike}
\bibliography{point}

\setcounter{equation}{0}
\renewcommand{\theequation}{A.\arabic{equation}}
\renewcommand{\thefigure}{A.\arabic{figure}}
\renewcommand{\thetable}{A.\arabic{figure}}

\section*{Appendix A: Proofs}

\subsection*{Proof of Lemma~\ref{lem.elliso}} 
Since both sides are zero-mean Gaussian processes, it will suffice to prove the equality of covariance structures (\citealp{adler:2007}, Ch.~1.2). 
But we have, for $\psi_1, \psi_2 \in L^2(B_\mu)$,
\begin{equation*}
\E[W_{B_\mu}(\psi_1)W_{B_\mu}(\psi_2)] = \langle \psi_1, \psi_2 \rangle_{B_\mu} = \langle\ell \psi_1, \ell \psi_2 \rangle_{B_\lambda}
= \E[W_{B_\lambda}(\ell\psi_1)W_{B_\lambda}(\ell\psi_2)],
\end{equation*}
which yields the desired result.

\subsection*{Proof of Lemma~\ref{lem.K}}
\begin{enumerate}
\item[(i)] We first prove that $R_{a,b}$ is bijective, that is, both injective and surjective.
For any $\varphi\in L^{2}(B_\lambda)$, we have
\allowdisplaybreaks{
\begin{align*}
R_{a,b}\circ R_{a,b}\, \varphi 
&= R_{a,b}\, \varphi -2\frac{\langle a-b,R_{a,b}\, \varphi\rangle_{B_\lambda}}{\|a-b\|^{2}_{B_\lambda}}(a-b)\\
&= \varphi -2\frac{\langle a-b,\varphi\rangle_{B_\lambda}}{\|a-b\|^{2}_{B_\lambda}}(a-b) -2\frac{\langle a-b,\varphi\rangle_{B_\lambda}}{\|a-b\|^{2}_{B_\lambda}}(a-b)\\ 
&\qquad+\ 2\frac{\langle a-b,2(a-b)\rangle_{B_\lambda}}{\|a-b\|^{2}_{B_\lambda}}\frac{\langle a-b,\varphi\rangle_{B_\lambda}}{\|a-b\|^{2}_{B_\lambda}}(a-b)\\
&=\varphi,
\end{align*}}%
which implies that $R_{a,b}$ is both injective, because $R_{a,b}\,\varphi_{1}=R_{a,b}\,\varphi_{2}$ implies $R_{a,b}\circ R_{a,b}\,\varphi_{1}=R_{a,b}\circ R_{a,b}\,\varphi_{2}$ hence $\varphi_{1}=\varphi_{2}$,
and surjective, because for all $\varphi\in L^{2}(B_\lambda)$, $R_{a,b}\,\varphi$ is mapped to $\varphi$.

Next, we prove that $R_{a,b}$ preserves the inner product.
We have that
\allowdisplaybreaks{
\begin{align*}
\lefteqn{\langle R_{a,b}\,\varphi_{1},R_{a,b}\,\varphi_{2}\rangle_{B_\lambda}} \qquad\\    
&=\langle \varphi_{1},\varphi_{2}\rangle_{B_\lambda} 
-\left\langle \varphi_{1},2\frac{\langle a-b,\varphi_{2}\rangle_{B_\lambda}}{\|a-b\|^{2}_{B_\lambda}}(a-b)\right\rangle_{B_\lambda}  \\ 
&\qquad -\left\langle \varphi_{2},2\frac{\langle a-b,\varphi_{1}\rangle_{B_\lambda}}{\|a-b\|^{2}_{B_\lambda}}(a-b)\right\rangle_{B_\lambda}  \\
&\qquad +\left\langle 2\frac{\langle a-b,\varphi_{1}\rangle_{B_\lambda}}{\|a-b\|^{2}_{B_\lambda}}(a-b),2\frac{\langle a-b,\varphi_{2}\rangle_{B_\lambda}}{\|a-b\|^{2}_{B_\lambda}}(a-b)\right\rangle_{B_\lambda}  \\
&=\langle \varphi_{1},\varphi_{2}\rangle_{B_\lambda}
-2\frac{\langle a-b,\varphi_{1}\rangle_{B_\lambda}\langle a-b,\varphi_{2}\rangle_{B_\lambda}}{\|a-b\|^{2}_{B_\lambda}}  \\ 
&\qquad-2\frac{\langle a-b,\varphi_{1}\rangle_{B_\lambda}\langle a-b,\varphi_{2}\rangle_{B_\lambda}}{\|a-b\|^{2}_{B_\lambda}}  
+4\frac{\langle a-b,\varphi_{1}\rangle_{B_\lambda}\langle a-b,\varphi_{2}\rangle_{B_\lambda}}{\|a-b\|^{2}_{B_\lambda}}  \\
&=\langle \varphi_{1},\varphi_{2}\rangle_{B_\lambda}.
\end{align*}}%
It follows that $R_{a,b}$ is both bijective and preserves the inner product, i.e., it is unitary.

\item[(ii)] Follows directly from the first part of the proof of~(i).

\item[(iii)] We have that
\begin{align*}
\langle R_{a,b}\,\varphi_{1},\varphi_{2}\rangle_{B_\lambda}    
&=\langle \varphi_{1},\varphi_{2}\rangle_{B_\lambda} 
-\left\langle 2\frac{\langle a-b,\varphi_{1}\rangle_{B_\lambda}}{\|a-b\|^{2}_{B_\lambda}}(a-b),\varphi_{2}\right\rangle_{B_\lambda}  \\ &=\langle \varphi_{1},\varphi_{2}\rangle_{B_\lambda}
-2\frac{\langle a-b,\varphi_{1}\rangle_{B_\lambda}\langle a-b,\varphi_{2}\rangle_{B_\lambda}}{\|a-b\|^{2}_{B_\lambda}}  \\
&=\langle \varphi_{1},\varphi_{2}\rangle_{B_\lambda}
-\left\langle \varphi_{1},2\frac{\langle a-b,\varphi_{2}\rangle_{B_\lambda}}{\|a-b\|^{2}_{B_\lambda}}(a-b)\right\rangle_{B_\lambda}  \\
&=\langle \varphi_{1},R_{a,b}\,\varphi_{2}\rangle_{B_\lambda}.
\end{align*}

\item[(iv)] Because both $a$ and $b$ have unit norm, one may verify that
\begin{align*}
R_{a,b}\, a 
&= a -2\frac{\langle a-b,a\rangle_{B_\lambda}}{\|a-b\|^{2}_{B_\lambda}}(a-b) = a -2\frac{\langle a-b,a\rangle_{B_\lambda}}{\langle a-b,a-b\rangle_{B_\lambda}}(a-b)\\
&= a -2\frac{1-\langle b,a\rangle_{B_\lambda}}{2-2\langle a,b\rangle_{B_\lambda}}(a-b) =b,
\end{align*}
and similarly for $R_{a,b}\, b$.
\end{enumerate}

\subsection*{Proof of Theorem~\ref{thm.main.1}}
As in Lemma~\ref{lem.elliso}, we are comparing two zero-mean Gaussian processes, and we will prove their equality in distribution by proving the equality of their covariance structures.

Let us denote by $\mathcal{L}(q_\mu)$ the linear span of $q_\mu$ in $L^2(B_\mu)$, that is, $\mathcal{L}(q_\mu) = \{c q_\mu : c \in \mathbb{R}\}$. Let us also denote by $\mathcal{L}^\perp(q_\mu)$ the orthogonal complement of $\mathcal{L}(q_\mu)$ in $L^2(B_\mu)$, that is, $\mathcal{L}^\perp(q_\mu) = \{\psi \in L^2(B_\mu) : \langle \psi, q_\mu \rangle_{B_\mu} = 0\}$. Then any $\psi \in L^2(B)$ can be uniquely decomposed as $\psi = \psi^{\|} + \psi^\perp$, with $\psi^{\|} \in \mathcal{L}(q_\mu)$ and $\psi^\perp \in \mathcal{L}^\perp(q_\mu)$. 

Now consider the process on the right-hand side of \eqref{eq.main.1}. We can write, for any $\psi \in L^2(B)$,
\begin{equation}\label{1.1}
V_{B_\mu}(\psi) = V_{B_\mu}(\psi^\|) + V_{B_\mu}(\psi^\perp).
\end{equation}
Since $\psi^\| = cq_\mu$ for some $c \in \mathbb{R}$, we can compute the first term on the right-hand side of \eqref{1.1} as
\[
V_{B_\mu}(\psi^\|) = cV_{B_\mu}(q_\mu) = c\big[W_{B_\mu}(q_\mu) - \langle q_\mu, q_\mu\rangle_{B_\mu} W_{B_\mu}(q_\mu) \big] = 0, 
\]
where we used the representation \eqref{f.i} for $V_{B_\mu}$, and the fact that $\langle q_\mu, q_\mu\rangle_{B_\mu}=1$; see \eqref{q.ortho}. As for the second term on the right-hand side of \eqref{1.1}, we have
\[
V_{B_\mu}(\psi^\perp) = W_{B_\mu}(\psi^\perp) - \langle \psi^\perp, q_\mu\rangle_{B_\mu} W_{B_\mu}(q_\mu) = W_{B_\mu}(\psi^\perp), 
\]
since $\langle \psi^\perp, q_\mu\rangle_{B_\mu} = 0$. It follows that $V_{B_\mu}(\psi) = W_{B_\mu}(\psi^\perp)$, and so the covariance structure of the process on the right-hand side of \eqref{eq.main.1} is given by
\begin{equation}\label{cov.1}
\E\Big[V_{B_\mu}(\psi_1)V_{B_\mu}(\psi_2)\Big] = \E\big[W_{B_\mu}(\psi_1^\perp)W_{B_\mu}(\psi_2^\perp)\big] = \langle \psi_1^\perp, \psi_2^\perp\rangle_{B_\mu}.
\end{equation}

Next, we consider the process on the left-hand side of \eqref{eq.main.1}. Once again, we have the decomposition
\begin{equation}\label{1.2}
V_{B_\lambda}(R_{q_\lambda,\ell q_\mu}\ell\psi) = V_{B_\lambda}(R_{q_\lambda,\ell q_\mu}\ell\psi^\|) + V_{B_\lambda}(R_{q_\lambda,\ell q_\mu}\ell\psi^\perp)
\end{equation}
for any $\psi \in L^2(B_\mu)$. Since $\psi^\| = cq_\mu$ for some $c \in \mathbb{R}$, we can compute the first term on the right-hand side of \eqref{1.2} as
\begin{align*}
\lefteqn{V_{B_\lambda}(R_{q_\lambda,\ell q_\mu}\ell\psi^\|)} \qquad \\ 
&= cV_{B_\lambda}(R_{q_\lambda,\ell q_\mu}\ell q_\mu) \stackrel{(*)}{=} cV_{B_\lambda}(q_\lambda) = c\big[W_{B_\lambda}(q_\lambda) - \langle q_\lambda, q_\lambda\rangle_{B_\lambda} W_{B_\lambda}(q_\lambda) \big] = 0,
\end{align*}
where the step $(*)$ follows from the fact that $R_{q_\lambda,\ell q_\mu}$ maps $\ell q_\mu$ to $q_\lambda$; see Lemma \ref{lem.K}(iv). As for the second term on the right-hand side of \eqref{1.2}, we have
\begin{align*}
\lefteqn{V_{B_\lambda}(R_{q_\lambda,\ell q_\mu}\ell\psi^\perp)} \qquad \\
&= W_{B_\mu}(R_{q_\lambda,\ell q_\mu}\ell\psi^\perp) - \langle R_{q_\lambda,\ell q_\mu}\ell\psi^\perp, q_\lambda\rangle_{B_\lambda} W_{B_\lambda}(q_\lambda) \stackrel{(\dagger)}{=} W_{B_\mu}(R_{q_\lambda,\ell q_\mu}\ell\psi^\perp),
\end{align*}
where the step $(\dagger)$ follows from the fact that $R_{q_\lambda,\ell q_\mu}\ell \, \cdot$ is a unitary operator, and therefore $R_{q_\lambda,\ell q_\mu}\ell\psi^\perp$ is perpendicular to $R_{q_\lambda,\ell q_\mu}\ell q_\mu = q_\lambda$.

It now follows that 
\[
V_{B_\lambda}(R_{q_\lambda,\ell q_\mu}\ell\psi) = W_{B_\lambda}(R_{q_\lambda,\ell q_\mu}\ell\psi^\perp),
\]
and so the covariance structure of the process on the right-hand side of \eqref{eq.main.1} is given by
\begin{equation}\label{cov.2}
\begin{split}
\E \Big[V_{B_\lambda}(R_{q_\lambda,\ell q_\mu}\ell\psi_1)\,V_{B_\lambda}(R_{q_\lambda,\ell q_\mu}\ell\psi_2) \Big] &= \E \big[ W_{B_\lambda}(R_{q_\lambda,\ell q_\mu}\ell\psi_1^\perp) W_{B_\lambda}(R_{q_\lambda,\ell q_\mu}\ell\psi_2^\perp) \big]\\
&= \langle R_{q_\lambda,\ell q_\mu}\ell\psi_1^\perp, R_{q_\lambda,\ell q_\mu}\ell\psi_2^\perp\rangle_{B_\lambda}\\
&= \langle \psi_1^\perp, \psi_2^\perp\rangle_{B_\mu},
\end{split}
\end{equation}
where the final equality follows from the fact that $R_{q_\lambda,\ell q_\mu}\ell \, \cdot$ is a unitary operator.

Comparing the final expressions in \eqref{cov.1} and \eqref{cov.2}, we see that the two processes in \eqref{eq.main.1} indeed have identical covariance structure, which completes the proof.

\subsection*{Proof of Theorem~\ref{thm.main.2}}
The proof is analogous to the proof of Theorem \ref{thm.main.1}, so we do not repeat it here in full. 
Once again, we start by defining $\mathcal{L}(q_\mu)$ to be the linear span of $q_\mu$ in $L^2(B_\mu)$, that is, $\mathcal{L}(q_\mu) = \{c^\top q_\mu : c \in \mathbb{R}^m \}$, and we define $\mathcal{L}^\perp(q_\mu)$ to be the orthogonal complement of $\mathcal{L}(q_\mu)$ in $L^2(B_\mu)$. 
Then we follow the proof of Theorem \ref{thm.main.1} almost \textit{verbatim}, adjusting for the fact that some scalars in that proof become vectors or matrices in the current setup, which does not change the arguments in the proof. 
We also use the fact that the unitary operator $U^{(m)}$ satisfies $U^{(m)} \ell q_\mu = q_\lambda$ by construction.

\subsection*{Proof of Theorem~\ref{thm.main}}
By Theorem \ref{thm:1} and Skorohod's representation theorem (\citealp{billingsley:1999}, Thm.~6.7), there exists a probability space supporting versions of $\widehat{W}_n$ and $V_{B_\lambda}$ satisfying
\begin{align}\label{Delta.W}
\sup_{t \in [0,T]} \big| \widehat{W}_n(t) - V_{B_\lambda}(t)\big| \to 0 \quad \text{a.s.}
\end{align}
We will show that on this space,
\begin{equation}\label{to.prove}
\sup_{t \in [0,S]} \bigg| \int_{0}^{T} \widehat{U}^{(m)} \widehat{\ell}\mathds{1}_{[0,t]}(s) \,\mathrm{d}\widehat{W}_n(s) - \int_{0}^{T}U^{(m)}\ell\mathds{1}_{[0,t]}(s)\,\mathrm{d}V_{B_\lambda}(s) \bigg| \stackrel{P}{\rightarrow} 0,
\end{equation}
which proves the theorem. 

Let us denote
\[
\Delta^V_n(s) := \widehat{W}_n(s) - V_{B_\lambda}(s), \quad s \in [0,T]. 
\]
Then, the difference of integrals on the left-hand side of \eqref{to.prove} can be written as
\[
\int_{0}^{T} \Delta^U_n(s,t) \, \dd \widehat{W}_n(s) + \int_{0}^{T} U^{(m)}\ell\mathds{1}_{[0,t]}(s) \, \dd \Delta^V_n(s),
\]
so the convergence in \eqref{to.prove} will follow from
\allowdisplaybreaks{
\begin{align}
\sup_{t \in [0,S]} \bigg| \int_{0}^{T} \Delta^U_n(s,t) \, \dd \widehat{W}_n(s) \bigg| \stackrel{P}{\rightarrow} 0, \label{to.prove1}\\
\sup_{t \in [0,S]} \bigg| \int_{0}^{T} U^{(m)}\ell\mathds{1}_{[0,t]}(s) \, \dd \Delta^V_n(s) \bigg| \stackrel{P}{\rightarrow} 0.\label{to.prove2}
\end{align}}\

Let us first consider the convergence \eqref{to.prove1}. 
Using integration by parts, we can write, for each $t \in [0,S]$,
\[
\int_{0}^{T} \Delta^U_n(s,t) \, \dd \widehat{W}_n(s) = \Delta^U_n(T,t) \widehat{W}_n(T) - \int_{0}^{T} \widehat{W}_n(s) \, \dd\Delta^U_n(s,t),
\]
and therefore the supremum in \eqref{to.prove1} is bounded by
\[
\sup_{t \in [0,T]} |\widehat{W}_n(t)| \cdot \Big( \sup_{t \in [0,S]} |\Delta_n^U(T,t)| + \sup_{t \in [0,S]} \mathrm{V}_0^T(\Delta^U_n(\cdot,t)) \Big),
\]
which converges to zero in probability, by Assumptions A4 (iii)--(iv), together with \eqref{Delta.W}. 

Similarly, the left-hand side of \eqref{to.prove2}, after using integration by parts, can be bounded by
\[
\sup_{t \in [0,T]} |\Delta^V_n(t)| \cdot \Big( \sup_{t \in [0,S]} |U^{(m)}\ell\mathds{1}_{[0,t]}(T)| + \sup_{t \in [0,S]} \mathrm{V}_0^T(U^{(m)}\ell\mathds{1}_{[0,t]}) \Big),
\]
which converges to zero in probability, by Assumptions A4 (i)--(ii), together with \eqref{Delta.W}.

\section*{Appendix B: Discussion of Assumption~A4}

In the case $m=1$, we note the following about conditions (i)--(iv) of Assumption A4.
\begin{enumerate}
\item[(i)] When $m=1$, this statement simply follows from the continuity of the mapping $t \mapsto U^{(1)}\ell\mathds{1}_{[0,t]}(T)$ over $t \in [0,S]$. Indeed we have, for $t \in [0,S]$,
\begin{align*}
U^{(1)}\ell\mathds{1}_{[0,t]}(T) &= R_{q_\lambda, \ell q_\mu}\ell\mathds{1}_{[0,t]}(T)\\
&= - 2 \,\frac{\int_0^t [q_\lambda(s) - \ell(s) q_\mu(s)] \ell(s) \beta(s)\, \dd s}{\int_0^T [q_\lambda(s) - \ell(s) q_\mu(s)]^2 \beta(s) \, \dd s}\,(q_\lambda(T) - \ell(T) q_\mu(T)),
\end{align*}
which is continuous in $t$.

\item[(ii)] For each $t \in [0,S]$, we can use the definition of $U^{(1)}\ell\mathds{1}_{[0,t]}$ and elementary properties of the total variation to write
\[
\mathrm{V}_0^T(U^{(1)}\ell\mathds{1}_{[0,t]})
\le \mathrm{V}_0^T(\ell) + 2 \,\frac{\int_0^T |q_\lambda(s) - \ell(s) q_\mu(s)| \ell(s) \beta(s)\, \dd s}{\int_0^T [q_\lambda(s) - \ell(s) q_\mu(s)]^2 \beta(s) \, \dd s}\mathrm{V}_0^T(q_\lambda - \ell q_\mu).
\]
So if the functions $\ell$, $q_\lambda$ and $q_\mu$ have derivatives that are bounded in absolute value over $[0,T]$, then the total variations on the right-hand side will be finite, and the statement A4(ii) will hold.

\item[(iii)] For $(s,t) \in [0,T]^2$, define $U^{(1)}_{\theta}\ell_{\theta}\mathds{1}_{[0,t]}(s)$ similar to $\widehat{U}^{(1)}\widehat{\ell}\mathds{1}_{[0,t]}(s)$, with the estimate $\widehat{\theta}$ replaced by arbitrary $\theta \in \Theta$. Suppose that the partial derivative
\[
\frac{\partial}{\partial \theta} \, U^{(1)}_{\theta}\ell_{\theta}\mathds{1}_{[0,t]}(T)
\]
exists and is bounded over $(t, \theta) \in [0,S] \times \Theta_0$. Now, since $\widehat\theta \in \Theta_0$ with probability tending to 1, there exists $\theta^* = \theta^*(t) \in \Theta_0$ with probability tending to 1 such that
\[
\Delta^U_n(T,t) = (\widehat\theta - \theta_0) \frac{\partial}{\partial \theta} \, U^{(1)}_{\theta}\ell_{\theta}\mathds{1}_{[0,t]}(T) \bigg|_{\theta=\theta^*},
\]
by the Mean Value Theorem. It follows that, with probability tending to 1,
\[
\sup_{t \in [0,S]} |\Delta^U_n(T,t)| \le (\widehat\theta - \theta_0) \sup_{(t, \theta) \in [0,S] \times \Theta_0} \bigg| \frac{\partial}{\partial \theta} \, U^{(1)}_{\theta}\ell_{\theta}\mathds{1}_{[0,t]}(T) \bigg|,
\]
where the difference $(\widehat\theta - \theta_0)$ is $o_P(1)$ under Assumptions A1--A3, and the supremum term is a finite constant. 
This will imply statement A4(iii).

\item[(iv)] Let $U^{(1)}_{\theta}\ell_{\theta}\mathds{1}_{[0,t]}(s)$ be as defined in part (iii) above. For any $(t,\theta) \in [0,T] \times \Theta_0$, the mapping $s \mapsto U^{(1)}_{\theta}\ell_{\theta}\mathds{1}_{[0,t]}(s)$ has a jump of size $\ell_\theta(t)$ at $s=t$, but is differentiable with respect to $s$ for $s \in [0,t) \cup (t,T]$, assuming $\alpha_\theta(s)$ and $\beta_\theta(s)$ of \eqref{eq:alpha-beta} are differentiable. Then, if we let
\[
D U^{(1)}_{\theta}\ell_{\theta}\mathds{1}_{[0,t]}(s) := \frac{\partial}{\partial s}U^{(1)}_{\theta}\ell_{\theta}\mathds{1}_{[0,t]}(s), 
\]
for $s \in [0,t) \cup (t,T]$, we can write
\begin{align*}
\sup_{t \in [0,S]} \mathrm{V}_0^T(\Delta^U_n(\cdot,t)) &\le T\sup_{t \in [0,S]}\sup_{s \in [0,t) \cup (t,T]} \big|D U^{(1)}_{\widehat\theta}\ell_{\widehat\theta}\mathds{1}_{[0,t]}(s) - D U^{(1)}_{\theta_0}\ell_{\theta_0}\mathds{1}_{[0,t]}(s)\big| \\
& \qquad {} + \sup_{t \in [0,S]} \big|\ell_{\widehat\theta}(t) - \ell_{\theta_0}(t)\big|.
\end{align*}
Now, if the functions $D U^{(1)}_{\theta}\ell_{\theta}\mathds{1}_{[0,t]}(s)$ and $\ell_\theta(s)$ have partial derivatives with respect to $\theta$ bounded over appropriate domains, then both of the supremum terms on the right-hand side will be $o_P(1)$, by a Mean Value Theorem argument as in (iii), and therefore A4(iv) will hold.
\end{enumerate}

\setcounter{equation}{0}
\renewcommand{\theequation}{C.\arabic{equation}}

\section*{Appendix C: Maximum Likelihood Estimation}

\subsection*{Aalen-Type Survival Processes}
We estimate the parameters of the Aalen model \eqref{eq:Aalen} with Weibull force of mortality $\varsigma_\theta$ as specified in \eqref{weibull.mort}, by maximizing the log-likelihood function \eqref{MLE}. This leads to the first-order conditions
\begin{equation}\label{MLE.aalen1}
\theta_1 = \bigg( \frac{1}{N_n(T)} \bigg[\sum_{i=1}^{N_n(T)} (t_0 + T_i)^{\theta_2} + (n-N_n(T))(t_0+T)^{\theta_2} -nt_0^{\theta_2} \bigg]\bigg)^{1/\theta_2}
\end{equation}
and
\begin{equation}\label{MLE.aalen2}
\begin{split}
\frac{N_n(T)}{\theta_2} + n \bigg( \frac{t_0}{\theta_1}\bigg)^{\theta_2}\log \bigg(\frac{t_0}{\theta_1}\bigg) + \sum_{i=1}^{N_n(T)} \bigg[ 1 - \bigg(\frac{t_0+T_i}{\theta_1} \bigg)^{\theta_2}\bigg]\log\bigg(\frac{t_0+T_i}{\theta_1} \bigg) \qquad \\ {} - (n-N_n(T))\bigg(\frac{t_0+T}{\theta_1} \bigg)^{\theta_2}\log\bigg(\frac{t_0+T}{\theta_1} \bigg) = 0.
\end{split}
\end{equation}
These equations can be solved numerically, by replacing $\theta_1$ in \eqref{MLE.aalen2} by the expression in \eqref{MLE.aalen1}, and then applying a standard root finding algorithm (such as \texttt{uniroot()} in \texttt{R}) to solve for $\theta_2$.

\subsection*{Software Reliability Models}
We estimate the parameters of the Jelinski-Moranda model \eqref{JM} by maximizing the log-likelihood function \eqref{MLE}, which leads to the first-order conditions
\begin{equation*}
\theta = \frac{N_n(T)}{\sum_{i=1}^{N_n(T)} T_i + (np-N_n(T))T},\qquad
0 = \sum_{i=1}^{N_n(T)} \frac{1}{np - i + 1} - \theta T,
\end{equation*}
where $T_1, \ldots, T_{N_n(T)}$ are the observed failure times in $[0,T]$. 
These equations can be solved numerically, by plugging in the first expression for $\theta$ into the second equation and then solving for $p$ via a standard root-finding algorithm.

In the case of the Littlewood model \eqref{Lw}, the log-likelihood function \eqref{MLE} takes the form
\begin{equation*}
\begin{split}
L(\theta_1, \theta_2, p) &= N_n(T) \log \theta_1 + \sum_{i=1}^{N_n(T)} \log(n p - i + 1) - \Big( 1 + \frac{\theta_1}{\theta_2} \Big) \sum_{i=1}^{N_n(T)} \log(1+ \theta_2 T_i)\\
&\qquad{} - \frac{\theta_1}{\theta_2}(np - N_n(T)) \log(1+\theta_2 T),    
\end{split}
\end{equation*}
where $T_1, \ldots, T_{N_n(T)}$ denote the observed failure times, as before. 
We maximize this function numerically (using \texttt{optim()} in \texttt{R}) to compute the maximum likelihood estimates.

\subsection*{Power Studies}
We estimate the parameters of the Aalen model \eqref{eq:Aalen} with Gompertz force of mortality $\varsigma_\theta$ as specified in \eqref{gomp.mort}, by maximizing the log-likelihood function \eqref{MLE}. 
This leads to the first-order conditions
\begin{equation*}
\theta_1 = N_n(T) \bigg(\sum_{i=1}^{N_n(T)} e^{\theta_2(t_0 + T_i)} + (n-N_n(T))e^{\theta_2(t_0+T)} -ne^{\theta_2 t_0} \bigg)^{-1}
\end{equation*}
and
\begin{equation*}
\begin{split}
\frac{N_n(T)}{\theta_2} + \theta_1 t_0 e^{\theta_2 t_0} N_n(T) + \sum_{i=1}^{N_n(T)} (t_0 + T_i) - \theta_1\sum_{i=1}^{N_n(T)}(t_0+T_i)e^{\theta_2(t_0 + T_i)} \qquad \\ {} + (n-N_n(T)) \theta_1 e^{\theta_2t_0} (t_0 - (t_0+T)e^{\theta_2 T}) = 0.
\end{split}
\end{equation*}
We solve these equations numerically, similar to the Weibull case.

\setcounter{equation}{0}
\renewcommand{\theequation}{D.\arabic{equation}}
\renewcommand{\thefigure}{D.\arabic{figure}}
\renewcommand{\thetable}{D.\arabic{figure}}

\section*{Appendix D: Additional Monte Carlo Simulations}

\subsection*{Aalen-Type Survival Processes with Random Censoring}

A common feature in survival and duration analysis is that of right-censoring (\citealp{kalbfleisch:2002}, Ch.~3; \citealp{khmaladze:2013}, Ch.~9). 
In the context of the Aalen-type survival models described in Section \ref{sec:Aalen}, it is natural to imagine that some of the $n$ individuals in the original population emigrate or move away from the population at random times, so that we cannot record their exact lifetimes in our data, but only lower bounds on their lifetimes (namely the times of their emigration). 
These incomplete lifetime measurements are referred to as right-censored lifetimes. 

The feature of right-censoring can be incorporated into the Aalen model as follows. 
We assume, as before, that the $n$ individuals in the population have i.i.d.\ lifetimes $T_1, \ldots, T_n$ generated from a parametric lifetime distribution $F_\theta$ with corresponding force of mortality $\varsigma_\theta$. 
We also assume that the $n$ individuals have i.i.d.\ ``censoring times'' $C_1, \ldots, C_n$ associated with them, which can be thought of as the random times of emigration out of the population. 
Note that we allow the random variables $C_i$ to take the value $\infty$ with positive probability, so that some individuals in the population never emigrate. 
Instead of observing the full lifetimes $T_1, \ldots, T_n$ of the $n$ individuals, we now observe the censored lifetimes
\[
\widetilde{T}_i = \min\{T_i, C_i\}, \quad i=1,\ldots,n,
\]
together with the indicator functions 
\[
\delta_i = \mathds{1}_{\{T_i = \widetilde{T}_i\}} = \mathds{1}_{\{T_i \le C_i\}}, \quad i=1,\ldots,n.
\]
In other words, for each individual $i$ in our original population, we observe either the time of death $T_i$ or the time of emigration $C_i$, whichever comes first. 
And we know for each observation $\widetilde{T}_i$ whether it is a full lifetime observation ($\delta_i=1$) or a right-censored one ($\delta_i=0$).

The point process $N_n(t)$ counting the number of deaths up to (and including) time $t$ now takes the form
\begin{equation}\label{eq:Ncens}
N_n(t) = \sum_{i=1}^n \mathds{1}_{\{T_i \le t\}} \delta_i, \quad 0 \le t \le T,
\end{equation}
where $[0,T]$ is our period of observation as before. Note that $N_n(t)$ does not count emigrations, but only deaths. 
The conditional intensity of $N_n(t)$ will no longer take the form (\ref{eq:Aalen}), but it will take the modified form
\begin{equation}
\lambda_{n,\theta}(t)=\varsigma_{\theta}(t)[n - \widetilde{N}_n(t-)],\qquad 0\leq t\leq T,
\label{eq:Aalencens}
\end{equation}
with $\widetilde{N}_n(t) = \sum_{i=1}^n \mathds{1}_{\{\widetilde{T}_i \le t\}}$, the number of deaths \emph{and emigrations} observed by time $t$. 
Note that both in \eqref{eq:Aalen} and in \eqref{eq:Aalencens}, the conditional intensity $\lambda_{n,\theta}(t)$ is given by the individual force of mortality at time $t$, $\varsigma_\theta(t)$, multiplied by the number of individuals still alive and present in the population ``just before'' time $t$.

As in Section~\ref{sec:Aalen}, we illustrate our approach with Weibull lifetimes, now with censoring. 
More precisely, we consider a population of $n$ individuals as before, each $t_0$ years old at time $t=0$, but instead of observing the full remaining lifetimes $T_i$ for each individual, we observe the censored random variables $\widetilde{T}_i = \min\{T_i,C_i\}$ along with the indicators $\delta_i = \mathds{1}_{\{\widetilde{T}_{i}=T_{i}\}}=\mathds{1}_{\{T_i \le C_i\}}$, where the $T_i$'s are i.i.d.\ Weibull distributed excess lifetimes with cdf as in (\ref{weibull.cdf}), and the $C_i$'s are i.i.d.\ censoring variables independent of the $T_i$'s. 
We assume that $C_i \sim IX + (1-I)Y$, where $I$ is a Bernoulli random variable with $\mathbb{P}[I=1] = 1 - \mathbb{P}[I=0] = p$, $X \sim \text{Exp}(\gamma)$ independent of $I$, and $Y$ is simply equal to $\infty$ (and $0 \cdot \infty$ is interpreted as 0). 

One can imagine that each individual in the population makes an independent decision to ``emigrate'' or ``move away'' from the area of observation with probability $p$. 
If an individual decides to emigrate, it takes them an $\text{Exp}(\gamma)$ amount of time to act on the decision, so only some of the individuals who make the decision get to emigrate before they die. 
For these individuals, we only observe the time of emigration rather than the time of death.

For simulations, we take $n=1000$ and $t_0=50$ as in Section~\ref{sec:Aalen}, as well as $\theta_1 = 86$ and $\theta_2 = 9$ for a realistic lifetime distribution. 
For censoring, we take $p=0.4$ and $\gamma = 1/15$, so that about 40\% of the original population makes the decision to emigrate, and it takes 15 years on average to act on this decision. 
We construct 1000 paths from the point process
\begin{equation*}
N_n(t) = \sum_{i=1}^n \mathds{1}_{\{\widetilde{T}_i \le t\}}\delta_i, \quad 0 \le t \le T,
\end{equation*}
which keeps track of the number of observed deaths (but not emigrations) in the population up to, and including, time $t$. We take $T=50$ as before. 
The conditional intensity process $\lambda_{n,\theta}$ associated with $N_n$ is given in (\ref{eq:Aalencens}).

From each generated sample path $N_n$, we estimate the parameters $\theta_1$ and $\theta_2$ by maximizing the log-likelihood function \eqref{MLE}, which in this case leads to the first-order conditions
\begin{equation}\label{MLE.aalencens1}
\theta_1 = \bigg( \frac{1}{N_n(T)} \bigg[\sum_{i=1}^{\widetilde{N}_n(T)} (t_0 + \widetilde{T}_i)^{\theta_2} + (n-\widetilde{N}_n(T))(t_0+T)^{\theta_2} -nt_0^{\theta_2} \bigg]\bigg)^{1/\theta_2}
\end{equation}
and
\begin{equation}\label{MLE.aalencens2}
\begin{split}
\frac{N_n(T)}{\theta_2} + n \bigg( \frac{t_0}{\theta_1}\bigg)^{\theta_2}\log \bigg(\frac{t_0}{\theta_1}\bigg) + \sum_{i=1}^{\widetilde{N}_n(T)} \bigg[ \delta_i - \bigg(\frac{t_0+\widetilde{T}_i}{\theta_1} \bigg)^{\theta_2}\bigg]\log\bigg(\frac{t_0+\widetilde{T}_i}{\theta_1} \bigg) \qquad \\ {} - (n-\widetilde{N}_n(T))\bigg(\frac{t_0+T}{\theta_1} \bigg)^{\theta_2}\log\bigg(\frac{t_0+T}{\theta_1} \bigg) = 0.
\end{split}
\end{equation}
Note that the system of equations \eqref{MLE.aalencens1}--\eqref{MLE.aalencens2} is similar to \eqref{MLE.aalen1}--\eqref{MLE.aalen2}, but not identical, due to the presence of right-censoring. 
Solving the first-order conditions numerically, we obtain the maximum likelihood estimates $\widehat\theta_1$ and $\widehat\theta_2$. 
Then from each realized path $N_n$, we construct the process $\widehat{W}_n$ described in \eqref{eq:Wnparametric}, as well as the transformed process $\widehat{\mathscr{T}}\,\widehat{W}_n$ in \eqref{eq.trhat}, on the interval $[0,T]$. 
The functions $\alpha_\theta$ and $\beta_\theta$ of (\ref{eq:alpha-beta}) in this case take the form:
\[
\alpha_\theta(t) = \bigg(\!\! -\frac{\theta_2}{\theta_1}, \, \frac{1}{\theta_2} + \log \bigg( \frac{t_0 + t}{\theta_1}\bigg)\bigg)^\top, \qquad \quad \beta_\theta(t) = f_\theta(t)(1-G(t)),
\]
where $G(t)$ denotes the cdf of the censoring variables $C_i$, so in this case:
\[
G(t) = p[1-\exp(-\gamma t)], \quad 0<t<\infty.
\]
To implement the empirical transformation $\widehat{\mathscr{T}}$, we use the functions $\alpha_{\widehat\theta}$ and $\widehat{\beta}_{\widehat\theta} = \frac{1}{n} \lambda_{n,\widehat\theta}$. 
Note that, since the censoring distribution $G$ is typically not known, we cannot simply use $\beta_{\widehat\theta}$ as an estimate for $\beta_\theta$, as we did in Section \ref{sec:Aalen}.

\begin{figure}[t]
\centering
\epsfig{file=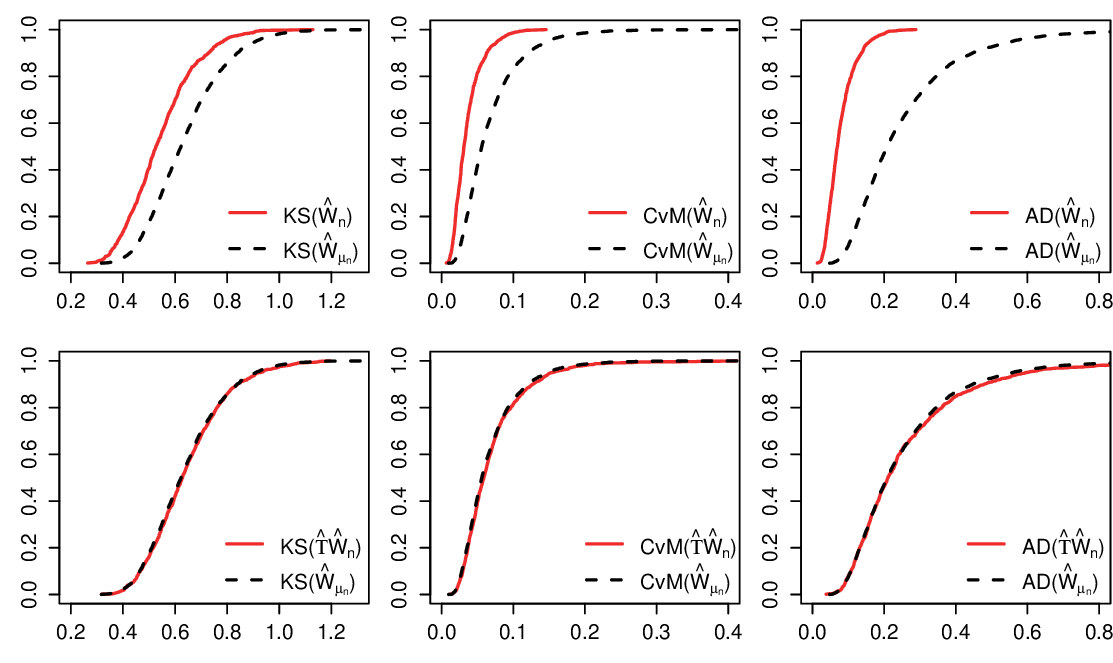,width=\textwidth}
\vspace{-15pt}
\caption{Testing for the Aalen model with censoring \eqref{eq:Aalencens} under the null hypothesis. The first row shows the empirical distribution functions of the three test statistics in \eqref{3stats}, as computed from $\widehat{W}_n$ and from $\widehat{W}_{\mu_n}$. The second row shows the empirical distribution functions of the same test statistics, as computed from $\widehat{\mathscr{T}}\widehat{W}_n$ and $\widehat{W}_{\mu_n}$. Unlike the empirical process $\widehat{W}_n$, the transformed process $\widehat{\mathscr{T}}\widehat{W}_n$ seems to behave identically to the target process $\widehat{W}_{\mu_n}$, in line with Theorem \ref{thm.main}.}
\label{fig.aalencens}
\end{figure}

From each generated sample path $N_n$, we construct the transformed process $\widehat{\mathscr{T}}\,\widehat{W}_n$ on $[0,50]$ as before, and we compare this process with $\widehat{W}_{\mu_n}$ using the three test statistics in \eqref{3stats}. 
Fig.~\ref{fig.aalencens} below is the ``censored case'' analog of Fig.~\ref{fig.aalen}: 
the first row compares the empirical cdfs of the three test statistics computed from the untransformed $\widehat{W}_n$ with the corresponding empirical cdfs computed from $\widehat{W}_{\mu_n}$, and the second row compares the empirical cdfs computed from the transformed $\widehat{\mathscr{T}}\,\widehat{W}_n$ with those computed from $\widehat{W}_{\mu_n}$. 
As expected, and in line with Theorem \ref{thm.main}, the empirical distributions are quite different before applying our transformation, but almost identical after. 

Table~\ref{table.AalenCens} is the analog of Table~\ref{table.Aalen}, and shows the observed rejection rates at 10\%, 5\% and 1\% significance levels among 1000 simulated samples. 
Note that the observed rejection rates are slightly higher than their expected values (albeit still close), most likely because censoring decreases the ``effective sample size'' and thereby increases the variances of all estimates. 

\begin{table}[b]
\begin{tabular}{c|cccc}
\textbf{Sig.\ level}   & \textbf{KS} & \textbf{CvM} & \textbf{AD} & \textbf{Expected} \\[2pt]
\hline\\[-7pt]
0.10 & 100 & 115 & 122 & 100\\[2pt]
\hline\\[-7pt]
0.05 & 56 & 60 & 64 & 50\\[2pt]
\hline\\[-7pt]
0.01 & 16 & 12 & 19 & 10\\[2pt]
\hline
\end{tabular}
\vspace{10pt}
\caption{Testing for the Aalen model with censoring \eqref{eq:Aalencens} under the null hypothesis. The table shows the number of times each test statistic lead to a rejection of the null hypothesis, among 1000 simulated paths. The observed rejection rates are close to the significance levels, in line with Theorem \ref{thm.main}.}\label{table.AalenCens}
\vspace{-10pt}
\end{table}

\subsection*{Mixture Cure Models}
Modified versions of the Aalen-type survival models discussed in Section \ref{sec:Aalen}, known as \emph{mixture cure models}, have found widespread use in medical studies, especially within cancer research, since the 1950s. 
The earliest mention of such models seems to appear in \cite{boag:1949};  we refer to \cite{maller:etal:2024} for a recent comprehensive overview.

As a simple example, consider a population of $n$ patients who have recently finished treatment for a certain type of cancer and have been declared cancer-free. 
While some of these patients are genuinely cured and will remain cancer-free in the future, others will experience a relapse of cancer some months or years after the treatment. 
If we are interested in modeling the length of the relapse period based on this population of $n$ individuals, we need to take into account the fact that the population consists of a mixture of ``uncured'' and ``cured'' individuals in unknown proportions. 
Among the uncured group, it is reasonable to model the relapse times as i.i.d.\ random variables with some common distribution $F$, while among the cured group, the relapse times will be effectively $\infty$. 
In particular, if we record the relapse times in this population during some observation period $[0,T]$, any patient who did not relapse by the end of the observation period might belong to either of the two groups: either they are cured and will never relapse, or they are uncured but their relapse time happens to be longer than the observation period.

To model this situation, we can imagine that an unknown proportion $p$ of the population remains uncured after the treatment and will relapse at random future times, while the remaining proportion $1-p$ is cured and will never relapse. 
We assume that among the uncured patients, the relapse times are independently generated from some absolutely continuous waiting time distribution $F_\theta$ (belonging to some pre-specified parametric family) with corresponding force of mortality $\varsigma_\theta$. 
Now, if $N_n(t)$ denotes the number of relapses observed in the population up to (and including) time $t$, then $N_n(t)$ will be a temporal point process with conditional intensity
\begin{equation}\label{eq:MixCure}
\lambda_{n,\theta,p}(t) = \varsigma_\theta(t)[np-N_n(t-)], \quad 0 \le t \le T,
\end{equation}
by arguments analogous to the ones given in Section~\ref{sec:Aalen}. 
Note that the form of \eqref{eq:MixCure} is identical to \eqref{eq:soft} in Section \ref{sec.software}, but the interpretation of the model components is different: here $\varsigma_\theta(t)$ denotes the force of mortality of the random relapse period instead of the waiting time until a failure, and $np$ denotes the initial number of uncured patients instead of the initial number of software faults. 
For known $p$, the model \eqref{eq:MixCure} reduces to \eqref{eq:Aalen}.

For simulations, we consider a population of $n=1000$ individuals with a ``cure rate'' of $1-p=0.25$. 
That is, we assume that 250 individuals in this population are cured and will never relapse, while the remaining 750 will relapse at random times. 
The relapse times are generated independently from a Weibull distribution with the following force of mortality and cdf:
\begin{equation}\label{weibull.cdf2}
\varsigma_\theta(t) = \frac{\theta_2}{\theta_1} \Big( \frac{t}{\theta_1} \Big)^{\theta_2 - 1}, \qquad F_\theta(t) = 1 - \exp\Big\{\!- \Big(\frac{t}{\theta_1}\Big)^{\theta_2}\Big\}, \qquad t \ge 0.
\end{equation}
We take $\theta_1 = 0.8$, $\theta_2 = 1.2$, and an observation period of $T=1$ year. These parameter choices lead to an average relapse time of about 0.75 years (or 9 months) among the uncured population, with a standard deviation of 0.63 years (or 7.5 months).

\begin{figure}[t]
\centering
\epsfig{file=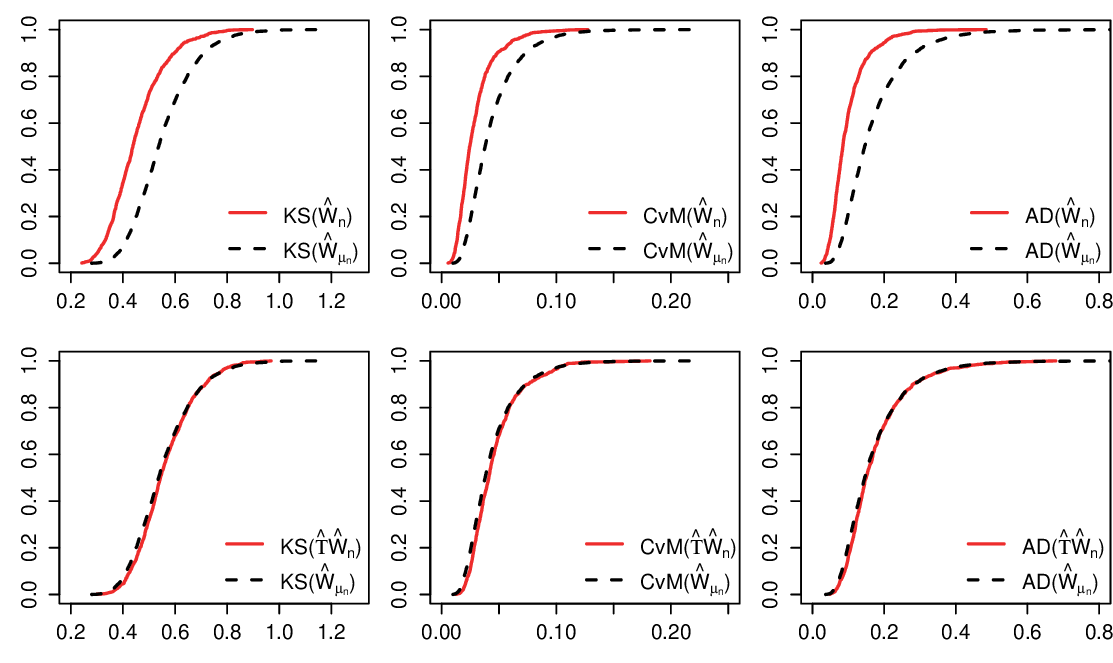,width=\textwidth}
\vspace{-15pt}
\caption{Testing for the mixture cure model \eqref{eq:MixCure} under the null hypothesis. The first row shows the empirical distribution functions of the three test statistics in \eqref{3stats}, as computed from $\widehat{W}_n$ and from $\widehat{W}_{\mu_n}$. The second row shows the empirical distribution functions of the same test statistics, as computed from $\widehat{\mathscr{T}}\widehat{W}_n$ and $\widehat{W}_{\mu_n}$. Unlike the empirical process $\widehat{W}_n$, the transformed process $\widehat{\mathscr{T}}\widehat{W}_n$ seems to behave identically to the target process $\widehat{W}_{\mu_n}$, in line with Theorem \ref{thm.main}.}
\label{fig.cure}
\end{figure}

The log-likelihood function \eqref{MLE} in this case takes the form
\begin{align*}
L(\theta_1, \theta_2 ,p) &= N_n(T) \log\Big(\frac{\theta_2}{\theta_1}\Big) - \sum_{i=1}^{N_n(T)} \Big( \frac{T_i}{\theta_1} \Big)^{\theta_2} + (\theta_2-1) \sum_{i=1}^{N_n(T)} \log\Big( \frac{T_i}{\theta_1} \Big) \\
& \qquad{} - (np-N_n(T))\Big( \frac{T}{\theta_1} \Big)^{\theta_2} + \sum_{i=1}^{N_n(T)} \log(np-i+1),
\end{align*}
where $T_1, \ldots, T_{N_n(T)}$ are the observed relapse times in $[0,T]$. 
Maximizing this log-likelihood numerically, we obtain the maximum likelihood estimates $\widehat\theta_1$, $\widehat\theta_2$ and $\widehat{p}$. 

Repeating the process above 1000 times, we obtain 1000 realized paths of $N_n(t)$ with corresponding sets of maximum likelihood estimates, and from each of these realized paths, we construct the processes $\widehat{W}_n$ and $\widehat{\mathscr{T}}\,\widehat{W}_n$ as in Sections~\ref{sec:Aalen} and \ref{sec.software}. 
The functions $\alpha_\theta$ and $\beta_\theta$ of (\ref{eq:alpha-beta}) take the form
\[
\alpha_{\theta,p}(t) = \bigg(\!\! -\frac{\theta_2}{\theta_1}, \, \frac{1}{\theta_2} + \log \Big( \frac{t}{\theta_1}\Big), \frac{1}{p (1-F_\theta(t))}\bigg)^\top, \qquad \quad \beta_{\theta,p}(t) = pf_\theta(t),
\]
and for the empirical transformation $\widehat{\mathscr{T}}$, we use $\alpha_{\widehat\theta, \widehat{p}}$ and $\beta_{\widehat\theta, \widehat{p}}$.

In Fig.~\ref{fig.cure}, we compare the empirical cdfs of the test statistics in (\ref{3stats}), computed from the untransformed $\widehat{W}_n$ as well as the transformed $\widehat{\mathscr{T}}\,\widehat{W}_n$, with the corresponding empirical cdfs computed from $\widehat{W}_{\mu_n}$ as before. As expected, the empirical cdfs of the test statistics in Fig.~\ref{fig.cure} are visibly different before the transformation (first row), but nearly identical after (second row). Table~\ref{table.Cure} is the analogue of Tables~\ref{table.Aalen}-\ref{table.JMLw}, and shows the observed rejection rates at 10\%, 5\% and 1\% significance levels among 1000 simulated samples. The observed rejection rates are once again close to their expected values.

\begin{table}[t]
\begin{tabular}{c|cccc}
\textbf{Sig.\ level}   & \textbf{KS} & \textbf{CvM} & \textbf{AD} & \textbf{Expected} \\[2pt]
\hline\\[-7pt]
0.10 & 106 & 103 & 96 & 100\\[2pt]
\hline\\[-7pt]
0.05 & 48 & 62 & 52 & 50\\[2pt]
\hline\\[-7pt]
0.01 & 5 & 8 & 12 & 10\\[2pt]
\hline
\end{tabular}
\vspace{10pt}
\caption{Testing for the mixture cure model \eqref{eq:MixCure} under the null hypothesis. The table shows the number of times each test statistic lead to a rejection of the null hypothesis, among 1000 simulated paths. The observed rejection rates are close to the significance levels, in line with Theorem \ref{thm.main}.}\label{table.Cure}
\vspace{-10pt}
\end{table}

\subsection*{Additional Power Studies}
\begin{figure}[t]
\centering
\epsfig{file=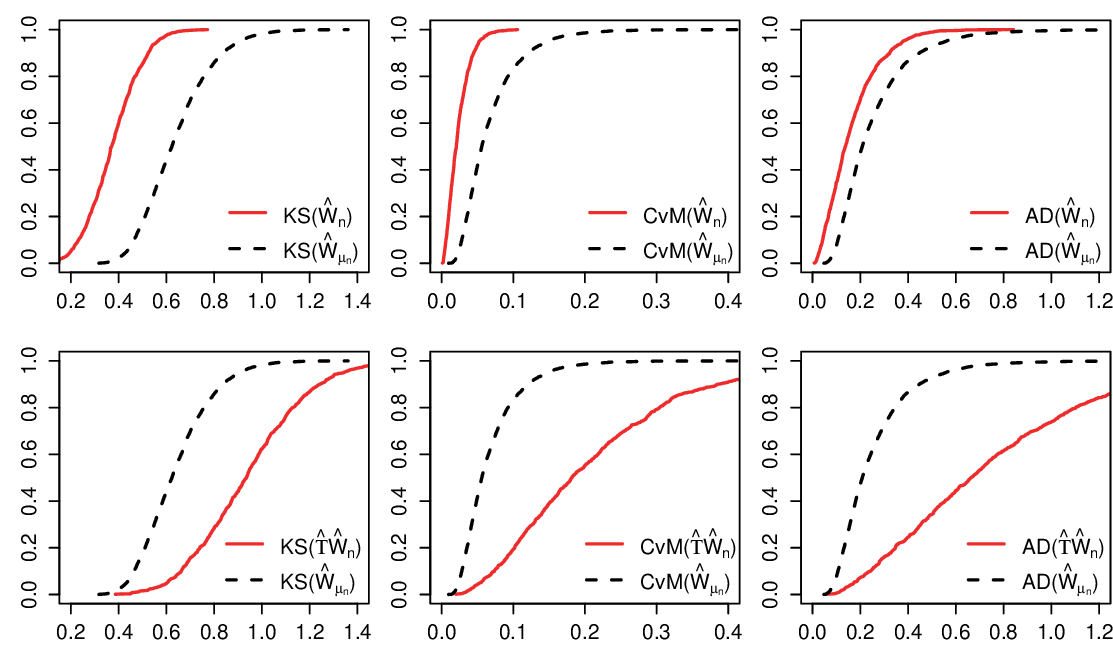,width=\textwidth}
\vspace{-15pt}
\caption{Testing for the Jelinski-Moranda model \eqref{JM} under the alternative hypothesis. The first row shows the empirical distribution functions of the three test statistics in \eqref{3stats}, as computed from $\widehat{W}_n$ and from $\widehat{W}_{\mu_n}$. The second row shows the empirical distribution functions of the same test statistics, as computed from $\widehat{\mathscr{T}}\widehat{W}_n$ and $\widehat{W}_{\mu_n}$. In contrast to the situation under the null hypothesis, the transformed process $\widehat{\mathscr{T}}\widehat{W}_n$ seems to be further away from the target process $\widehat{W}_{\mu_n}$ than the un-transformed $\widehat{W}_n$, implying good power properties.}
\label{fig.power2}
\end{figure}

In the case of software reliability models, we generate 1000 independent sample paths from the process $N_n(t)$ with the Littlewood conditional intensity (\ref{Lw}), with parameters chosen as in Section~\ref{sec.software}: $n=10{,}000$, $\theta_1=4$, $\theta_2=1$, $p=0.1$. The length of the observation period is $T=1$. From each generated sample, we construct the processes $\widehat{W}_n$ and $\widehat{\mathscr{T}}\,\widehat{W}_n$, now under the false assumption of a Jelinski-Moranda conditional intensity (\ref{JM}). 
Then we compare these two processes with the target process $\widehat{W}_{\mu_n}$ using the three test statistics in (\ref{3stats}), as before. 
See Fig.~\ref{fig.power2} for the results. 
Once again, we observe that the transformed process $\widehat{\mathscr{T}}\,\widehat{W}_n$ is not any closer to the target process $\widehat{W}_{\mu_n}$ than the un-transformed process $\widehat{W}_n$. Table~\ref{table.power2} displays the observed rejection rates (out of 1000) at 10\%, 5\% and 1\% significance levels. As suggested by Fig.\ \ref{fig.power2}, the rejection rates are much higher than what would be expected under the null hypothesis.

\begin{table}[t]
\begin{tabular}{c|cccc}
\textbf{Sig.\ level}   & \textbf{KS} & \textbf{CvM} & \textbf{AD} & \begin{tabular}[c]{@{}c@{}}\textbf{Expected}\\\textbf{(under null)}\end{tabular} \\[2pt]
\hline\\[-7pt]
0.10 & 656 & 719 & 707 & 100\\[2pt]
\hline\\[-7pt]
0.05 & 542 & 619 & 596 & 50\\[2pt]
\hline\\[-7pt]
0.01 & 306 & 407 & 370 & 10\\[2pt]
\hline
\end{tabular}
\vspace{10pt}
\caption{Testing for the Jelinski-Moranda model \eqref{JM} under the alternative hypothesis. The table shows the number of times each test statistic lead to a rejection of the null hypothesis, among 1000 simulated paths. In contrast to the situation under the null hypothesis, the observed rejection rates are much higher than the corresponding significance levels, implying good power properties.}\label{table.power2}
\vspace{-10pt}
\end{table}

\end{document}